\documentclass[12pt,english]{amsart}
\usepackage[T1]{fontenc}
\usepackage[latin9]{inputenc}
\usepackage{verbatim}
\usepackage{amsmath}
\usepackage{amssymb}
\usepackage{esint}

\makeatletter

\usepackage{amsfonts}
\usepackage{amscd}
\usepackage{amsfonts}
\providecommand{\U}[1]{\protect\rule{.1in}{.1in}}
\providecommand{\U}[1]{\protect\rule{.1in}{.1in}}

\newtheorem{theorem}{Theorem}[section]
\newtheorem{lemma}[theorem]{Lemma}
\newtheorem{proposition}[theorem]{Proposition}

\newtheorem{example}[theorem]{Example}
\newtheorem{definition}[theorem]{Definition}
\newtheorem{corollary}[theorem]{Corollary}
\newtheorem{remark}[theorem]{Remark}
\newtheorem{notation}[theorem]{Notation}
\newcommand{\be}{\begin{equation}}
\newcommand{\ee}{\end{equation}}

\newcommand{\norm}[1]{\Vert #1\Vert}

\newcommand{\cB}{\mathcal{B}}

\newcommand{\cD}{\mathcal{D}}

\newcommand{\cF}{\mathcal{F}}

\makeatother

\usepackage{babel}

\begin{document}

\title[Noncommutative function theory]{Progress in noncommutative function
theory}

\author{Paul S. Muhly}
\address{Department of Mathematics\\
 University of Iowa\\
 Iowa City, IA 52242}
\email{pmuhly@math.uiowa.edu}
\thanks{The research of the first author was supported by a grant from
the U.S.-Israel Binational Science
Foundation.}

\author{Baruch Solel}
\address{Department of Mathematics\\
 Technion\\
 32000 Haifa, Israel}
\email{mabaruch@techunix.technion.ac.il}
\thanks{The research of the second author was supported by the
U.S.-Israel Binational Science Foundation
and by the Lowengart Research Fund.}

\date{}
\dedicatory{To our esteemed friend and teacher,\\Richard V. Kadison,\\ on the
happy occasion of his $85^{th}$ birthday}
\begin{abstract}
In this expository paper we describe the study of certain non-self-adjoint
operator algebras, the Hardy algebras, and their representation theory.
We view these algebras as algebras of (operator valued) functions
on their spaces of representations. We will show that these spaces
of representations can be parameterized as unit balls of certain $W^{*}$-correspondences
and the functions can be viewed as Schur class operator functions
on these balls. We will provide evidence to show that the elements
in these (non commutative) Hardy algebras behave very much like bounded
analytic functions and the study of these algebras should be viewed
as noncommutative function theory.
\end{abstract}
\maketitle

\section{Introduction}

In this paper we shall introduce the tensor and Hardy operator algebras
and discuss how to study them as algebras of operator valued functions
on their representation spaces.

Tensor algebras associated with a bimodule over a ring have been studied
extensively in a purely algebraic setting. This class of algebras
has proved to be very important. In fact, every finite dimensional
algebra is a quotient of a tensor algebra.

Looking for a similar class of operator algebras, we were led by the
pioneering work of Pimsner \cite{Pi97} to study operator algebras
associated with $C^{*}$-correspondences. A $C^{*}$-correspondence
is, roughly, a bimodule over a $C^{*}$-algebra $M$ that is also
a (right) Hilbert $C^{*}$-module (see Section 2 below for more details).

These operator algebras, which we call tensor algebras, are subalgebras
of the $C^{*}$-algebras studied by Pimsner and are closely related
to them. Both are generated (as a norm-closed algebra and as a $C^{*}$-algebra,
respectively) by {}``shifts\textquotedbl{} on the Fock space of the
correspondence. In fact, the Cuntz-Pimsner algebra associated to a
given correspondence $E$ can be shown to be the {}``minimal\textquotedbl{}
$C^{*}$-algebra that is generated by the tensor algebra of $E$.
We shall not need this here but details can be found in \cite{MS98}
and \cite{KK06}.

In this paper we shall take the $C^{*}$-algebra $M$ to be a $W^{*}$-algebra
and assume that $E$ is a $W^{*}$-correspondence (details and definitions
are in the next section). This allows us to take the ultra-weak closure
of the tensor algebra. We call this ultra-weakly closed algebra the
Hardy algebra associated with the correspondence. As we shall see
below, the Hardy algebra that we get in the simplest case (where $M=\mathbb{C}=E$)
is simply the classical Hardy algebra $H^{\infty}(\mathbb{D})$. The
Hardy algebras associated with general $W^{*}$-correspondences are
the main object of our study here.

When studying the representations of the tensor algebras, we realized
that they can be parameterized by points in the closed unit balls
of certain $W^{*}$-correspondences. This fact will be exploited when
we view the elements of the tensor or the Hardy algebras as functions
on the representation space.

Considering the elements of an algebra as functions on the set of
its representations is not new, of course. It was done in a purely
algebraic setting and in the setting of Banach or $C^{*}$-algebras.
But, as we shall see, the fact that the representation space here
can be viewed as a unit ball of a $W^{*}$-correspondence, will allow
us to view these algebras as generalizations of algebras of holomorphic
functions on the disc $\mathbb{D}$ in $\mathbb{C}$.

In the next section we shall define the tensor and Hardy algebras
and describe some of their basic properties. As we shall see, these
algebras are generated by a copy of the $W^{*}$-algebra $M$ and
a copy of the correspondence $E$.

In Section 3 we study the representation theory of the tensor and
Hardy algebras. We shall first discuss the (completely contractive)
representations of the tensor algebras. For this, we fix a normal
representation $\sigma$ of $M$ on a Hilbert space $H$ and then
show that all the representations of the algebra whose restriction
to the copy of $M$ is $\sigma$ can be parameterized by the points
of the closed unit ball of a certain $W^{*}$-correspondence (that
we call the $\sigma$-dual of $E$ and write $E^{\sigma}$ for it).

In order to study the (completely contractive, ultra-weakly continuous)
representations of the Hardy algebra, we have to find out what representations
of the tensor algebra can be extended to such representations of the
Hardy algebra. This is done in Subsection 3.2. In this way, we identify
the ultra-weakly continuous representations of the Hardy algebra $H^{\infty}(E)$
as a subset of the closed unit ball of $E^{\sigma}$ that contains
the open unit ball. We write $AC(E^{\sigma})$ for this set.

Given a point $\eta\in AC(E^{\sigma})$, the associated representation
of the Hardy algebra will be denoted $\eta^{*}\times\sigma$ and,
given an element $X\in H^{\infty}(E)$, we write \begin{equation}
\widehat{X}(\eta^{*})=(\eta^{*}\times\sigma)(X).\end{equation}
 The reason for evaluating the function at $\eta^{*}$ and not at
$\eta$ is technical and will be clarified later.

We, thus, obtain the transform $X\mapsto\widehat{X}$, where $\widehat{X}$
is an operator valued function. We have already discussed the domain
of these functions. In Section 4 we discuss the nature of these functions
and we shall see that, up to a constant multiple, they form a natural
generalization of the classical Schur class functions. We present
two characterizations of these functions and call them Schur class
operator functions.

In the last two sections we take a closer look at the transform $X\mapsto\widehat{X}$
. In Section 5 we discuss the kernel of the transform and in the last
section we note that we are really dealing with several transforms:
for each normal representation $\sigma$ of $M$ we get a different
transform and we discuss the relationships among them.

Along the way, we present several results that demonstrate our main
point of view: These Hardy algebras form a useful analogue of the
algebra of holomorphic functions on the disc $\mathbb{D}$ and their
study can be seen as noncommutative function theory.

\section{Introducing the tensor and the Hardy algebras}

Before we introduce the algebras, we describe the setup. Throughout
this paper, $M$ will denote a fixed $W^{*}$-algebra. We do not preclude
the possibility that $M$ may be finite dimensional. Indeed, the situation
when $M=\mathbb{C}^{d}$ can be very interesting (even for $d=1$).
However, we want to think of $M$ abstractly, as a $C^{*}$-algebra
that is a dual space, without regard to any Hilbert space on which
$M$ might be represented. The weak-$*$ topology on a $W^{*}$-algebra
or on any of its weak-$*$ closed subspaces will be referred to as
the \emph{ultra-weak} topology.

To eliminate unnecessary technicalities, we shall always assume $M$
is $\sigma$-finite in the sense that every family of mutually orthogonal
projections in $M$ is countable. Alternatively, to say $M$ is $\sigma$-finite
is to say that $M$ has a faithful normal representation on a separable
Hilbert space. So, unless explicitly indicated otherwise, every Hilbert
space we consider will be assumed to be separable.

In addition, $E$ will denote a $W^{*}$-correspondence over $M$
in the sense of \cite{MSHardy}. For the definition, recall first
that a (right) Hilbert $C^{*}$-module over $M$ is a right module
$E$ over $M$ that is also equipped with an $M$-valued inner product.
More explicitly, we have a function $\langle\cdot,\cdot\rangle:E\times E\rightarrow M$
such that, for $\xi,\eta\in E$ and $a\in M$,
\begin{enumerate}
\item [1.] $\zeta\mapsto\langle\xi,\zeta\rangle$ is linear,
\item [2.] $\langle\xi,\eta a\rangle=\langle\xi,\eta\rangle a$
\item [3.] $\langle\xi,\eta\rangle$=$\langle\eta,\xi\rangle^{*}$,
\item [4.] $\langle\xi,\xi\rangle\geq0$, with $\langle\xi,\xi\rangle=0$
only if $\xi=0$, and
\item [5.] $E$ is complete in the norm $\norm{\xi}:=\norm{\langle\xi,\xi\rangle}^{1/2}$.
\end{enumerate}
Such a $C^{*}$-module is said to be self-dual provided each (right)
module map $\Phi$ from $E$ into $M$ is induced by a vector in $E$,
i.e., there is an $\eta\in E$ such that $\Phi(\xi)=\langle\eta,\xi\rangle$,
for all $\xi\in E$.

A self-dual Hilbert $C^{*}$-module $E$ over a $W^{*}$-algebra $M$
is said to be a $W^{*}$-module. Our basic reference for Hilbert $C^{*}$-
and $W^{*}$- modules is \cite{MT05}. It is shown in \cite[Proposition 3.3.4]{MT05}
that when $E$ is a self-dual Hilbert module over a $W^{*}$-algebra
$M$, then $E$ must be a dual space. In fact, it may be viewed as
an ultra-weakly closed subspace of a $W^{*}$-algebra. Further, every
continuous module map on $E$ is adjointable \cite[Corollary 3.3.2]{MT05}
and the algebra $\mathcal{L}(E)$ consisting of all continuous module
maps on $E$ is a $W^{*}$-algebra \cite[Proposition 3.3.4]{MT05}.

Given a $W^{*}$-module $E$ over $M$ and a normal $^{*}$-representation
$\sigma$ of $M$ on a Hilbert space $H$, one can define on the algebraic
tensor product, $E\otimes H$, a (scalar valued) inner product that
satisfies $\langle\xi\otimes h,\eta\otimes k\rangle=\langle h,\sigma(\langle\xi,\eta\rangle_{E})k\rangle_{H}$.
The completion of this inner-product space is a Hilbert space and
we write $E\otimes_{\sigma}H$ for it. One can then define the induced
representation $\sigma^{E}$ of $\mathcal{L}(E)$ on $E\otimes_{\sigma}H$
by \begin{equation}
\sigma^{E}(X)(\xi\otimes h)=X\xi\otimes h\;,\quad X\in\mathcal{L}(E),\;\xi\in E,\; h\in H.\end{equation}
 We shall also write $X\otimes I_{H}$ for $\sigma^{E}(X)$.

\begin{definition} Let $E$ be a $W^{*}$-module over the $W^{*}$-algebra
$M$. We say that $E$ is a $W^{*}$-correspondence over $M$ if there
is an ultra-weakly continuous $*$-representation $\varphi:M\to\mathcal{L}(E)$
such that $E$ becomes a bimodule over $M$ where the left action
of $M$ is determined by $\varphi_{E}$ (or simply $\varphi$), $a\cdot\xi=\varphi(a)\xi$.

We shall assume that $E$ is \emph{essential} or \emph{non-degenerate}
as a left $M$-module. This is the same as assuming that $\varphi$
is unital.

We also shall assume that our $W^{*}$-correspondences are countably
generated as self-dual Hilbert modules over their coefficient algebras.
This is equivalent to assuming that $\mathcal{L}(E)$ is $\sigma$-finite.
\end{definition}

\begin{example}

(Basic Example) If $M=\mathbb{C}$, then a $W^{*}$-correspondence
over $M$ is simply a Hilbert space. \end{example}

\begin{example}\label{quivercorrespondence}
Let $G=(G^{0},G^{1},r,s)$ be a directed graph. For simplicity we
assume that $G$ is finite. Thus both the set of vertices, $G^{0}$,
and the set of edges, $G^{1}$, are finite; and $r,s:G^{1}\rightarrow G^{0}$
are the range and source maps. We set $M=\ell^{\infty}(G^{0})$ (so
that $M$ is simply $\mathbb{C}^{n}$, for some $n$, viewed as a
$W^{*}$-algebra), and we set $E=\ell^{\infty}(G^{1})$. Then we endow
$E$ with the structure of a $W^{*}$-correspondence via the formulas:
\[
(\varphi(a)\xi b)(e)=a(r(e))\xi(e)b(s(e))\;,\quad a,b\in M,\;\;\xi\in E\;,\; e\in G^{1},\]
and \[
\langle\xi,\eta\rangle(v)=\sum_{s(e)=v}\langle\xi(e),\eta(e)\rangle\;,\quad\xi,\eta\in E\;,\;\; v\in G^{0}.\]
 One can easily check that every $W^{*}$-correspondence over a finite
dimensional commutative $W^{*}$-algebra is associated in this way
with a finite directed graph. \end{example} \begin{example} Let
$M$ be an arbitrary ($\sigma$-finite) $W^{*}$-algebra and let $\alpha:M\rightarrow M$
be a normal $^{*}$-endomorphism. Let $E=M$ (as a vector space) with
right action given by multiplication, left action given by $\varphi=\alpha$
and inner product $\langle\xi,\eta\rangle:=\xi^{*}\eta$. We denote
this correspondence by $_{\alpha}M$. (If $\alpha$ is the identity,
we write simply $M$ for this correspondence). \end{example} \begin{example}
Let $\Phi$ be a normal, contractive, completely positive map on the
$W^{*}$-algebra $M$. Write $E=M\otimes_{\Phi}M$. This is the $W^{*}$-correspondence
obtained as the self-dual completion of the algebraic tensor product
$M\otimes M$ with the inner product defined by $\langle a\otimes b,c\otimes d\rangle=b^{*}\Phi(a^{*}c)d$
and the bimodule structure defined by left and right multiplication:
$\varphi(c)(a\otimes b)d=ca\otimes bd$. This correspondence was used
by Popa \cite{P86}, Mingo \cite{Mi89}, Anantharam-Delarouche \cite{AnD90}
and others to study the map $\Phi$. It is referred to as the GNS
correspondence of $\Phi$. If $\Phi$ is an automorphism, $M\otimes_{\Phi}M$
is isomorphic to $_{\Phi}M$.

\end{example}

Along with $E$, we may form the ($W^{*}$-)tensor powers of $E$,
$E^{\otimes n}$. They will be understood to be the self-dual completions
of the $C^{*}$-tensor powers of $E$. Recall that the $C^{*}$-tensor
product of two correspondences $E$ and $F$ over $M$ is the completion
of the algebraic (balanced) tensor product $E\otimes F$ with respect
to the inner product \[
\langle\xi_{1}\otimes\zeta_{1},\xi_{2}\otimes\zeta_{2}\rangle=\langle\zeta_{1},\varphi_{F}(\langle\xi_{1},\xi_{2}\rangle_{E})\zeta_{2}\rangle_{F}\;,\quad\xi_{1},\xi_{2}\in E,\;\;\zeta_{1},\zeta_{2}\in F\]
 and the bimodule structure is defined by \[
\varphi_{E\otimes F}(a)(\xi\otimes\zeta)b=\varphi_{E}(a)\xi\otimes\zeta b,\;\quad\xi\in E,\;\zeta\in F,\; a,b\in M.\]

Likewise, the Fock space over $E$, $\mathcal{F}(E)$, will be the
self-dual completion of the Hilbert $C^{*}$-module direct sum of
the $E^{\otimes n}$: \[
\mathcal{F}(E)=M\oplus E\oplus E^{\otimes2}\oplus E^{\otimes3}\oplus\cdots\]
 We view $\mathcal{F}(E)$ as a $W^{*}$-correspondence over $M$,
where the left and right actions of $M$ are the obvious ones, i.e.,
the diagonal actions, and we shall write $\varphi_{\infty}$ for the
left diagonal action of $M$. Thus, for $\xi_{1}\otimes\xi_{2}\otimes\cdots\otimes\xi_{k}\in E^{\otimes k}$
and $a\in M$, \[
\varphi_{\infty}(a)(\xi_{1}\otimes\xi_{2}\otimes\cdots\otimes\xi_{k})=(\varphi(a)\xi_{1})\otimes\xi_{2}\otimes\cdots\otimes\xi_{k}.\]

For $\xi\in E$, we shall write $T_{\xi}$ for the so-called \emph{creation
operator} on $\mathcal{F}(E)$ defined by the formula $T_{\xi}\eta=\xi\otimes\eta$,
$\eta\in\mathcal{F}(E)$. It is easy to see that $T_{\xi}$ is in
$\mathcal{L}(\mathcal{F}(E))$ with norm $\Vert\xi\Vert$, and that
$T_{\xi}^{*}$ annihilates $M$, as a summand of $\mathcal{F}(E)$,
while on elements of the form $\zeta\otimes\eta$, $\zeta\in E,$ $\eta\in\mathcal{F}(E)$,
it is given by the formula \[
T_{\xi}^{*}(\zeta\otimes\eta):=\varphi_{\infty}(\langle\xi,\zeta\rangle)\eta.\]

We are now ready to define the operator algebras.

\begin{definition}\label{def:Tensor and Hardy Algs}If $E$ is a
$W^{*}$-correspondence over a $W^{*}$-algebra $M$, then \textbf{the
tensor algebra} of $E$, denoted $\mathcal{T}_{+}(E)$, is defined
to be the norm-closed subalgebra of $\mathcal{L}(\mathcal{F}(E))$
generated by $\varphi_{\infty}(M)$ and $\{T_{\xi}\mid\xi\in E\}$.
The \textbf{Hardy algebra} of $E$, denoted $H^{\infty}(E)$, is defined to
be the ultra-weak closure in $\mathcal{L}(\mathcal{F}(E))$ of $\mathcal{T}_{+}(E)$.
\end{definition}

\begin{example} If $M=E=\mathbb{C}$, the Fock correspondence is
the Hilbert space $\ell^{2}$ and, for $\xi=1\in\mathbb{C}$, $T_{1}$
is the unilateral shift. The tensor algebra in this case is the norm-closed
algebra generated by the shift and can be identified with the disc
algebra $A(\mathbb{D})$. The Hardy algebra is its $w^{*}$-closure
and can be identified with $H^{\infty}(\mathbb{D})$. \end{example}
It will be useful to bear this example in mind as we proceed because
our algebras, in general, can be viewed as noncommutative analogues
of the disc and the (classical) Hardy algebras.

\begin{example} If $M=\mathbb{C}$ and $E=\mathbb{C}^{d}$, then
the Fock correspondence is the Hilbert space $\ell^{2}(\mathbb{F}_{d}^{+})$
where $\mathbb{F}_{d}^{+}$ is the free semigroup on $d$ generators.
Letting $\{e_{i}:1\leq i\leq d\}$ be the standard orthonormal basis
of $E=\mathbb{C}^{d}$, we see that the tensor algebra is generated
(as a norm-closed algebra) by the $d$ shifts $\{T_{e_{i}}:1\leq i\leq d\}$
and the Hardy algebra is its $w^{*}$-closure. These algebras were
studied extensively by Popescu (e.g. \cite{gP96}), Davidson and Pitts
(e.g. \cite{DP98}) and others. Popescu denoted this tensor algebra
$\mathcal{A}_{d}$ (and called it the noncommutative disc algebra).
The Hardy algebra was denoted $F_{d}^{\infty}$ by Popescu and $\mathcal{L}_{d}$
by Davidson and Pitts. \end{example}

More examples are given in \cite{MSHardy} and discussed in detail
there.

An important tool used in the analysis of $\mathcal{T}_{+}(E)$ and
$H^{\infty}(E)$ is the {}``spectral theory of the gauge automorphism
group''. What we need is developed in detail in \cite[Section 2]{MSHardy}.
Here we merely recall the essentials. The reader should keep in mind
that its primary role is to allow us to handle in an analytic way
the natural gradings that the Fock space and the Hardy algebra have.
Let $P_{n}$ denote the projection of $\mathcal{F}(E)$ onto $E^{\otimes n}$.
Then $P_{n}\in\mathcal{L}(\mathcal{F}(E))$ and the series\[
W_{t}:=\sum_{n=0}^{\infty}e^{int}P_{n}\]
 converges in the ultra-weak topology on $\mathcal{L}(\mathcal{F}(E))$.
The family $\{W_{t}\}_{t\in\mathbb{R}}$ is an ultra-weakly continuous,
$2\pi$-periodic unitary representation of $\mathbb{R}$ in $\mathcal{L}(\mathcal{F}(E))$.
Further, if $\{\gamma_{t}\}_{t\in\mathbb{R}}$ is defined by the formula
$\gamma_{t}=Ad(W_{t})$, then $\{\gamma_{t}\}_{t\in\mathbb{R}}$ is
an ultra-weakly continuous group of $*$-automorphisms of $\mathcal{L}(\mathcal{F}(E))$
that leaves invariant $\mathcal{T}_{+}(E)$ and $H^{\infty}(E)$.
Indeed, the subalgebra of $H^{\infty}(E)$ fixed by $\{\gamma_{t}\}_{t\in\mathbb{R}}$
is $\varphi_{\infty}(M)$ and $\gamma_{t}(T_{\xi})=e^{-it}T_{\xi}$,
$\xi\in E$. Associated with $\{\gamma_{t}\}_{t\in\mathbb{R}}$ we
have the {}``Fourier coefficient operators'' $\{\Phi_{j}\}_{j\in\mathbb{Z}}$
on $\mathcal{L}(\mathcal{F}(E))$, which are defined by the formula

\begin{equation}
\Phi_{j}(a):=\frac{1}{2\pi}\int_{0}^{2\pi}e^{-int}\gamma_{t}(a)\, dt,\qquad a\in\mathcal{L}(\mathcal{F}(E)),\label{FourierOperators}\end{equation}
 where the integral converges in the ultra-weak topology. An alternate
formula for $\Phi_{j}$ is \[
\Phi_{j}(a)=\sum_{k\in\mathbb{Z}}P_{k+j}aP_{k}.\]
 Each $\Phi_{j}$ leaves $H^{\infty}(E)$ invariant and, in particular,
$\Phi_{j}(T_{\xi_{1}}T_{\xi_{2}}\cdots T_{\xi_{n}})=T_{\xi_{1}}T_{\xi_{2}}\cdots T_{\xi_{n}}$
if and only if $n=j$ and zero otherwise. Associated with the $\Phi_{j}$
are the {}``arithmetic mean operators'' $\{\Sigma_{k}\}_{k\geq1}$
that are defined by the formula\[
\Sigma_{k}(a):=\sum_{|j|<k}(1-\frac{|j|}{k})\Phi_{j}(a),\]
 $a\in\mathcal{L}(\mathcal{F}(E))$. For $a\in\mathcal{L}(\mathcal{F}(E))$,
$\lim_{k\to\infty}\Sigma_{k}(a)=a$, where the limit is taken in the
ultra-weak topology.

Note that, for $X\in H^{\infty}(E)$ and $k\geq1$, $\Phi_{k}(X)=T_{\xi_{k}}$
for some $\xi_{k}\in E^{\otimes k}$ and $\Phi_{0}(X)=\varphi_{\infty}(a)$
for some $a\in M$. We can write the {}``Fourier expansion\textquotedbl{}
of $X$ \begin{equation}
X\sim\Phi_{0}(X)+\Phi_{1}(X)+\Phi_{2}(X)+\cdots=\varphi_{\infty}(a)+T_{\xi_{1}}+T_{\xi_{2}}+\cdots.\label{eq:Taylor series.}\end{equation}

\section{The representations of the tensor and the Hardy algebras}

We now turn to describe the representation theory of $\mathcal{T}_{+}(E)$
and $H^{\infty}(E)$ . Details for what we describe are presented
in Section 2 of \cite{MSHardy} and in \cite{MSAbsCont}.

\subsection{Representations of the tensor algebras}

We start by discussing the representations of the tensor algebra $\mathcal{T}_{+}(E)$.

We shall consider only completely contractive representations and,
in fact, only those completely contractive representations of $\mathcal{T}_{+}(E)$
with the property that $\rho\circ\varphi_{\infty}$ is an ultra-weakly
continuous representation of $M$. This is not a significant restriction.
In particular, it is not a restriction at all, if $H$ is assumed
to be separable, since every $C^{*}$-representation of a $\sigma$-finite
$W^{*}$-algebra on a separable Hilbert space is automatically ultra-weakly
continuous \cite[Theorem V.5.1]{mT79}.

Note that, in the purely algebraic setting, where $M$ is a ring and
$E$ is an $M$-bimodule, the representations of the (algebraic) tensor
algebra are given by bimodule maps on $E$.

Here, suppose $\rho$ is a completely contractive representation of
$\mathcal{T}_{+}(E)$ on a Hilbert space $H$ as above, then $\sigma:=\rho\circ\varphi_{\infty}$
is a normal $^{*}$-representation of $M$ on $H$ and $\rho$ defines
a bimodule map $T$ from $E$ to $B(H)$ by the formula\[
T(\xi):=\rho(T_{\xi}).\]
 To say that $T(\cdot)$ is a bimodule map means simply that $T(\varphi(a)\xi b)=\sigma(a)T(\xi)\sigma(b)$
for all $a,b\in M$ and for all $\xi\in E$. The assumption that $\rho$
is completely contractive guarantees that $T$ is completely contractive
with respect to the unique operator space structure on $E$ that arises
from viewing $E$ as a corner of its linking algebra.

\begin{definition} \label{Definition1.12}Let $E$ be a $W^{\ast}$-correspondence
over a $W^{*}$-algebra $M$. Then:
\begin{enumerate}
\item A \emph{completely contractive covariant representation }of $E$ on
a Hilbert space $H$ is a pair $(T,\sigma)$, where

\begin{enumerate}
\item $\sigma$ is a normal $\ast$-representation of $N$ in $B(H)$.
\item $T$ is a linear, completely contractive map from $E$ to $B(H)$
that is continuous in the $\sigma$-topology of \cite{BDH88} on $E$
and the ultraweak topology on $B(H).$
\item $T$ is a bimodule map in the sense that
\[T(\varphi(a)\xi b)=\sigma(a)T(\xi)\sigma(b),\quad
\xi\in E, \mbox{ } a,b\in M.\]
\end{enumerate}
\item A completely contractive covariant representation $(T,\sigma)$ of
$E$ in $B(H)$ is called \emph{isometric }in case \begin{equation}
T(\xi)^{\ast}T(\eta)=\sigma(\langle\xi,\eta\rangle)\label{isometric}\end{equation}
 for all $\xi,\eta\in E$.
\end{enumerate}
\end{definition}


The discussion above shows that every completely contractive representation
$\rho$ of $\mathcal{T}_{+}(E)$ on $H$ gives rise to a completely
contractive covariant representation of $E$ on $H$. The following theorem
shows that the converse also holds and it can be viewed as a generalized
von Neumann inequality.

\begin{theorem} \label{Theorem310MS98}Let $E$ be a $W^{\ast}$-correspondence
over a von Neumann algebra $M$. To every completely contractive covariant
representation, $(T,\sigma)$, of $E$ there is a unique completely
contractive representation $\rho$ of the tensor algebra $\mathcal{T}_{+}(E)$
that satisfies \[
\rho(T_{\xi})=T(\xi)\;\;\;\xi\in E\]
 and\[
\rho(\varphi_{\infty}(a))=\sigma(a)\;\;\; a\in M.\]
 The map $(T,\sigma)\mapsto\rho$ is a bijection between the set of
all completely contractive covariant representations of $E$ and all
completely contractive (algebra) representations of $\mathcal{T}_{+}(E)$
whose restrictions to $\varphi_{\infty}(M)$ are continuous with respect
to the ultraweak topology on $\mathcal{L}(\mathcal{F}(E))$. \end{theorem}
\begin{definition} \label{integratedform}If $(T,\sigma)$ is a completely
contractive covariant representation of a $W^{\ast}$-correspondence
$E$ over a von Neumann algebra $M$, we call the representation $\rho$
of $\mathcal{T}_{+}(E)$ described in Theorem \ref{Theorem310MS98}
the \emph{integrated form} of $(T,\sigma)$ and write $\rho=T\times\sigma$.
\end{definition}

As we showed in \cite[ Lemmas 3.4--3.6]{MS98}, and in \cite{MSHardy},
if a completely contractive covariant representation, $(T,\sigma)$,
of $E$ in $B(H)$ is given, then it determines a contraction $\tilde{T}:E\otimes_{\sigma}H\rightarrow H$
defined by the formula $\tilde{T}(\eta\otimes h):=T(\eta)h$, $\eta\otimes h\in E\otimes_{\sigma}H$.
The operator $\tilde{T}$ intertwines the representation $\sigma$
on $H$ and the induced representation $\sigma^{E}:=\varphi(\cdot)\otimes I_{H}$
of $M$ on $E\otimes_{\sigma}H$; i.e. \begin{equation}
\tilde{T}(\varphi(\cdot)\otimes I)=\sigma(\cdot)\tilde{T}.\label{covariance}\end{equation}
 In fact we have the following lemma from \cite[Lemma 2.16]{MSHardy},.

\begin{lemma} \label{CovRep}The map $(T,\sigma)\rightarrow\tilde{T}$
is a bijection between all completely contractive covariant representations
$(T,\sigma)$ of $E$ on the Hilbert space $H$ and contractive operators
$\tilde{T}:E\otimes_{\sigma}H\rightarrow H$ that satisfy equation
(\ref{covariance}). Given such a $\tilde{T}$ satisfying this equation,
$T$, defined by the formula $T(\xi)h:=\tilde{T}(\xi\otimes h)$,
together with $\sigma$ is a completely contractive covariant representation
of $E$ on $H$. Further, $(T,\sigma)$ is isometric if and only if
$\tilde{T}$ is an isometry. \end{lemma}

Associated with $(T,\sigma)$ we also have maps $\tilde{T}_{n}:E^{\otimes n}\otimes H\rightarrow H$
defined by $\tilde{T}_{n}(\xi_{1}\otimes\xi_{2}\cdots\otimes\xi_{n}\otimes h)=T(\xi_{1})T(\xi_{2})\cdots T(\xi_{n})h$.

Now fix a normal representation $\sigma$ of $M$ on a Hilbert space
$H$. The discussion above shows that the set of all the completely
contractive representations $\rho$ of the tensor algebra $\mathcal{T}_{+}(E)$
that satisfy $\rho\circ\varphi_{\infty}=\sigma$ (roughly speaking,
$\rho$, restricted to $M$ is $\sigma$) can be parameterized by
the contractions $\tilde{T}\in B(E\otimes_{\sigma}H,H)$ that satisfy
the intertwining relation (\ref{covariance}). This is, of course,
the same as saying that this set of representations are parameterized
by the adjoints $\tilde{T}^{*}$. The reason that we prefer to consider
the adjoints is that the set of all the maps $\tilde{T}^{*}$ satisfying
relation (\ref{covariance}) can be given the structure of a $W^{*}$-correspondence
as the following proposition shows.

\begin{proposition}\label{corres} Let $E$ be a $W^{*}$-correspondence
over the $W^{*}$-algebra $M$ and let $\sigma$ be a normal representation
of $M$ on the Hilbert space $H$. Write $E^{\sigma}$ for the space
of all bounded maps $\eta:H\rightarrow E\otimes_{\sigma}H$ that satisfy
\begin{equation}
\eta\sigma(a)=(\varphi_{\infty}(a)\otimes I_{H})\eta\;,\quad a\in M.\end{equation}

With respect to the action of $\sigma(M)^{\prime}$ and the $\sigma(M)^{\prime}$-valued
inner product defined as follows, $E^{\sigma}$ becomes a $W^{\ast}$-correspondence
over $\sigma(M)^{\prime}$: For $X,Y\in\sigma(M)^{\prime}$, and $T\in E^{\sigma}$,
$X\cdot T\cdot Y:=(I\otimes X)TY$, and for $T,S\in E^{\sigma}$,
$\langle T,S\rangle:=T^{\ast}S$. \end{proposition}

\begin{definition} The $W^{*}$-correspondence of Proposition~\ref{corres}
will be called the $\sigma$-dual of $E$. \end{definition}

From equation \eqref{covariance} we see that $\widetilde{T}^{*}$
lies in the space we have denoted $E^{\sigma}$. So, if we write $\mathbb{D}(E^{\sigma})$
for the open unit ball in $E^{\sigma}$ and $\overline{\mathbb{D}(E^{\sigma})}$
for its norm closure, then all the completely contractive representations
$\rho$ of $\mathcal{T}_{+}(E)$ such that $\rho\circ\varphi_{\infty}=\sigma$
are parametrized bijectively by $\overline{\mathbb{D}(E^{\sigma*})}=\overline{\mathbb{D}(E^{\sigma})^{*}}=\overline{\mathbb{D}(E^{\sigma})}^{*}$.

\begin{example} In the special case when $(E,M)$ is $(\mathbb{C}^{d},\mathbb{C}),$
a representation $\sigma$ of $\mathbb{C}$ on a Hilbert space $H$
is quite simple; it does the only thing it can: $\sigma(c)h=ch$,
$h\in H$, and $c\in\mathbb{C}$. In this setting, $E\otimes_{\sigma}H$
is just the direct sum of $d$ copies of $H$ and $\widetilde{T}$
is simply a $d$-tuple of operators $(T_{1},T_{2},\ldots,T_{d})$
such that $\Vert\sum_{i}T_{i}T_{i}^{*}\Vert\leq1$, i.e. $\widetilde{T}$
is a row contraction. The map $T$, then, is given by the formula
$T(\xi)=\sum\xi_{i}T_{i}$, where $\xi=(\xi_{1},\xi_{2},\cdots,\xi_{d})^{\top}\in\mathbb{C}^{d}$.
The space $E^{\sigma}$ is column space over $B(H)$, $\mathbf{C}_{d}(B(H))$,
and $\mathbb{D}(E^{\sigma})$ is simply the unit ball in $\mathbf{C}_{d}(B(H))$.
\end{example}


It follows from Pimsner's analysis that $(T,\sigma)$ is isometric
if and only if $T\times\sigma$ is the restriction to $\mathcal{T}_{+}(E)$
of a $C^{*}$-representation of the $C^{*}$-subalgebra $\mathcal{T}(E)$
of $\mathcal{L}(\mathcal{F}(E))$ generated by $\mathcal{T}_{+}(E)$.
This $C^{*}$-algebra is called the \emph{Toeplitz algebra} of $E$.

A special kind of isometric covariant representations that will play
an important role here are constructed as follows. Let $\pi_{0}:M\to B(H_{0})$
be a normal representation of $M$ on the Hilbert space $H_{0}$,
and let $H=\mathcal{F}(E)\otimes_{\pi_{0}}H_{0}$. Set $\sigma:=\pi^{\mathcal{F}(E)}\circ\varphi_{\infty}=\varphi_{\infty}(\cdot)\otimes I_{H_{0}}$,
and define $S:E\to B(H)$ by the formula $S(\xi)=T_{\xi}\otimes I_{H_{0}}$,
$\xi\in E$. Then it is immediate that $(S,\sigma)$ is an isometric
covariant representation and we say that it is \emph{induced by $\pi_{0}$}.
We also will say $S\times\sigma$ is induced by $\pi_{0}$. In fact,
\begin{equation}
S\times\sigma=\pi_{0}^{\mathcal{F}(E)}|\mathcal{T}_{+}(E).\label{induced}\end{equation}

In a sense that will become clear, an induced representations should
be viewed as a generalization of a unilateral shift where the representation
$\pi_{0}$ plays the role of the multiplicity of the shift.

An induced isometric covariant representation has the property that
$\widetilde{S_{n}}\widetilde{S_{n}^{*}}\rightarrow0$ strongly as
$n\to\infty$ because $\widetilde{S_{n}}\widetilde{S_{n}^{*}}$ is
the projection onto $\sum_{k\geq n}E^{\otimes k}\otimes_{\pi_{0}}H_{0}$.
In general, an isometric covariant representation $(S,\sigma)$ and
its integrated form are called \emph{pure} if $\widetilde{S_{n}}\widetilde{S_{n}^{*}}\to0$
strongly as $n\to\infty$.

Corollary 2.10 of \cite{MS99} shows that every pure isometric covariant
representation of $(E,M)$ is unitarily equivalent to an isometric
covariant representation that is induced by a normal representation
of $M$. We therefore will usually say simply that a pure isometric
covariant representation \emph{is} induced. In Theorem 2.9 of \cite{MS99}
we proved a generalization of the Wold decomposition theorem that
asserts that every isometric covariant representation of $(E,M)$
decomposes as the direct sum of an induced isometric covariant representation
of $(E,M)$ and an isometric representation of $(E,M)$ that is both
isometric and fully coisometric.

We will need an analogue of a unilateral shift of infinite multiplicity.
For that, we shall fix, once and for all, a representation $(S_{0},\sigma_{0})$
that is induced by a faithful normal representation $\pi$ of $M$
that has \emph{infinite multiplicity}. That is, $(S_{0},\sigma_{0})$
acts on a Hilbert space of the form $\cF(E)\otimes_{\pi}K_{0}$, where
$\pi:M\to B(K_{0})$ is an infinite ampliation of a faithful normal
representation of $M$. Then $\sigma_{0}:=\pi^{\cF(E)}\circ\varphi_{\infty}$,
while $S_{0}(\xi):=T_{\xi}\otimes I_{K_{0}}$, $\xi\in E$. The following
proposition shows the uniqueness and the special role of this representation.

\begin{proposition}\label{Lemma:Universal}The representation $(S_{0},\sigma_{0})$
is unique up to unitary equivalence and every induced isometric covariant
representation of $(E,M)$ is unitarily equivalent (in a natural way)
to a restriction of $(S_{0},\sigma_{0})$ to a subspace of the form
$\cF(E)\otimes_{\pi}\mathfrak{K}$, where $\mathfrak{K}$ is a subspace
of $K_{0}$ that reduces $\pi$.\end{proposition}

\begin{definition}\label{Definition:Universal}We shall refer to
$(S_{0},\sigma_{0})$ as \emph{the universal induced covariant representation}
of $(E,M)$. \end{definition}By Proposition \ref{Lemma:Universal},
$(S_{0},\sigma_{0})$ does not really depend on the choice of representation
$\pi$ used to define it. It will serve the purpose in our theory
that the unilateral shift of infinite multiplicity serves in the structure
theory of single operators on Hilbert space.

A key tool in our theory is the following result that we proved as
\cite[Theorem 2.8]{MSHardy}.

\begin{theorem}\label{isomdil} Let $(T,\sigma)$ be a completely
contractive covariant representation of $(E,M)$ on a Hilbert space
$H$. Then there is an isometric covariant representation $(V,\tau)$
of $(E,M)$ acting on a Hilbert space $K$ containing $H$ such that
if $P$ denotes the projection of $K$ onto $H$, then
\begin{enumerate}
\item $P$ commutes with $\tau(M)$ and $\tau(a)P=\sigma(a)P$, $a\in M$,
and
\item for all $\eta\in E$, $V(\eta)^{*}$ leaves $H$ invariant and $PV(\eta)P=T(\eta)P$.
\end{enumerate}
The representation $(V,\tau)$ may be chosen so that the smallest
subspace of $K$ that contains $H$ and is invariant under both $\tau(M)$
and $V(E)$, is all of $K$. When this is done, $(V,\tau)$ is unique
up to unitary equivalence and is called \textbf{\emph{the minimal
isometric dilation}}\emph{ of} $(T,\sigma)$. \end{theorem}

Note that, in the notation of the theorem, we have \begin{equation}
T\times\sigma=P(V\times\tau)P.\label{compression}\end{equation}
 Thus, the representation $T\times\sigma$ is a compression, onto
a coinvariant subspace, of the representation $V\times\tau$.

Another result that will be important when studying the representations
of the tensor and the Hardy algebras is our version of the commutant
lifting theorem. This theorem was proved in \cite{MS98} and can be
stated as follows (see \cite[Theorems 2.6 and 2.7]{MSAbsCont}).

\begin{theorem}\label{CLT}For $i=1,2$, let $(T_{i},\sigma_{i})$
be a completely contractive covariant representation of $(E,M)$ on
a Hilbert space $H_{i}$, let $(V_{i},\tau_{i})$ be the minimal isometric
dilation of $(T_{i},\sigma_{i})$ acting on the space $K_{i}$, and
let $P_{i}$ be the orthogonal projection of $K_{i}$ onto $H_{i}$.
Then, given an operator $X\in B(H_{1},H_{2})$ that intertwines the
representations $T_{1}\times\sigma_{1}$ and $T_{2}\times\sigma_{2}$,
there is an operator $Y\in B(K_{1},K_{2})$ such that
\begin{enumerate}
\item [(1)] $Y$ intertwines the representations $V_{1}\times\tau_{1}$
and $V_{2}\times\tau_{2}$,
\item [(2)] $X=P_{2}YP_{1}$,
\item [(3)] $YH_{1}^{\perp}\subseteq H_{2}^{\perp}$ and
\item [(4)] $\norm{Y}=\norm{X}$.
\end{enumerate}
\end{theorem}

We end this section with a discussion of the representations of the
tensor algebras associated with directed graphs (see Example~\ref{quivercorrespondence}).

\begin{example}\label{quiverrepresentations} Let $G$ and $E$ as
described in Example~\ref{quivercorrespondence}. Write $E(G)$ for
$E$. The algebra $H^{\infty}(E)$ in this case will be written $H^{\infty}(G)$.
In the literature, $H^{\infty}(G)$ is sometimes denoted $\mathcal{L}_{G}$.
It is the ultraweak closure of the tensor algebra $\mathcal{T}_{+}(E(G))$
acting on the Fock space $\mathcal{F}(E(G))$. For $e\in G^{1}$,
let $\delta_{e}$ be the $\delta$-function at $e$, i.e., $\delta_{e}(e^{\prime})=1$
if $e=e^{\prime}$ and is zero otherwise. Then $T_{\delta_{e}}$ is
a partial isometry that we denote by $S_{e}$. Also, for $v\in G^{0}$,
$P_{v}$ is defined to be $\varphi_{\infty}(\delta_{v})$. Then each
$P_{v}$ is a projection and it is an easy matter to see that the
families $\{S_{e}:e\in G^{1}\}$ and $\{P_{v}:v\in G^{0}\}$ form
a \emph{Cuntz-Toeplitz family} in the sense that the following conditions
are satisfied:
\begin{enumerate}
\item [(i)] $P_{v}P_{u}=0$ if $u\neq v$,
\item [(ii)] $S_{e}^{*}S_{f}=0$ if $e\neq f$
\item [(iii)] $S_{e}^{*}S_{e}=P_{s(e)}$ and
\item [(iv)] $\sum_{r(e)=v}S_{e}S_{e}^{*}\leq P_{v}$ for all $v\in G^{0}$.
\end{enumerate}
The algebra $\mathcal{T}_{+}(E(G))$ was first defined and studied
in \cite{pM97}, providing examples of the theory developed in \cite{MS98}.
It was called a quiver algebra there because in pure algebra, directed
graphs are called quivers. The properties of quiver algebras were
further developed in \cite{MS99}. In \cite{KP}, the focus was on
$H^{\infty}(G)$ and the authors called this algebra a free semigroupoid
algebras. Both algebras are often represented as algebras of operators
on $l_{2}(G^{\ast})$ (where $G^{*}$ is the set of all finite paths
in $G$), and it will be helpful to understand how this is done, from
the perspective of this note. Let $H_{0}$ be a Hilbert space whose
dimension equals the number of vertices, let $\{e_{v}|\; v\in G^{0}\}$
be a fixed orthonormal basis for $H_{0}$ and let $\pi_{0}$ be the
diagonal representation of $M=\ell^{\infty}(G^{0})$ on $H_{0}$.
Then $l_{2}(G^{\ast})$ is isomorphic to $\mathcal{F}(E(G))\otimes_{\pi_{0}}H_{0}$
where the isomorphism maps an element $\xi_{\alpha}$ of the standard
orthonormal basis of $l_{2}(G^{\ast})$ to $\delta_{\alpha}\otimes e_{s(e_k)}$
(where, for a finite path $\alpha=e_{1}\cdots e_{k}$, $\delta_{\alpha}=\delta_{e_{1}}\otimes\cdots\otimes\delta_{e_{k}}\in E^{\otimes k}$).
The partial isometries $S_{e}$ can then be viewed as the shift operators
$S_{e}\xi_{\alpha}=\xi_{e\alpha}$. Thus, the representations of $\mathcal{T}_{+}(E(G))$
and $H^{\infty}(G)$ on $l_{2}(G^{\ast})$ are just the representations
induced by $\pi_{0}$.

Quite generally, a completely contractive covariant representation
of $E(G)$ on a Hilbert space $H$ is given by a representation $\sigma$
of $M=\ell^{\infty}(G^{0})$ on $H$ and by a contractive map $\tilde{T}:E\otimes_{\sigma}H\rightarrow H$
satisfying equation (\ref{covariance}). The representation $\sigma$
is given by the projections $Q_{v}=\sigma(\delta_{v})$ whose sum
is $I$. Also, from $\tilde{T}$ we may define maps $T(e)\in B(H)$
by the equation $T(e)h=\tilde{T}(\delta_{e}\otimes h)$ and it is
easy to check that $\tilde{T}\tilde{T}^{\ast}=\sum_{e}T(e)T(e)^{\ast}$
and $T(e)=Q_{r(e)}T(e)Q_{s(e)}$. Thus to every completely contractive
representation of the quiver algebra $\mathcal{T}_{+}(E(G))$ we associate
a family $\{T(e)|e\in G^{1}\}$ of maps on $H$ that satisfy $\sum_{e}T(e)T(e)^{\ast}\leq I$
and $T(e)=Q_{r(e)}T(e)Q_{s(e)}$. Conversely, every such family defines
a representation, written $T\times\sigma$ (or $\tilde{T}\times\sigma$),
satisfying $(T\times\sigma)(S_{e})=T(e)$ and $(T\times\sigma)(P_{v})=Q_{v}$.

Now we fix $\sigma$ to be $\pi_{0}$ and write $H$ in place of $H_{0}$.
So that, in this case, each projection $Q_{v}$ is one dimensional
(with range equal to $\mathbb{C}e_{v}$). Then obviously $\sigma(M)^{\prime}=\sigma(M)$.
To describe the $\sigma$-dual of $E$, write $G^{-1}$ for the directed
graph obtained from $G$ by reversing all arrows, so that $s(e^{-1})=r(e)$
and $r(e^{-1})=s(e)$. Sometimes $G^{-1}$ is denoted $G^{op}$ and
is called the opposite graph. Note that the Hilbert space $E\otimes_{\sigma}H_{0}$
is spanned by the orthonormal basis $\{\delta_{e}\otimes e_{s(e)}\}$.
Fix $\eta\in E^{\sigma}$ and note that its covariance property implies
that, for every $e\in G^{1}$, $\eta^{\ast}(\delta_{e}\otimes e_{s(e)})=\eta^{\ast}(\delta_{r(e)}\delta_{e}\otimes e_{s(e)})=Q_{r(e)}\eta^{\ast}(\delta_{e}\otimes e_{s(e)})=\overline{\eta(e^{-1})}e_{r(e)}$
for some $\overline{\eta(e^{-1})}\in\mathbb{C}$. The reason for the
{}``strange\textquotedbl{} way of writing that scalar is that now
we can view $\eta$ as an element of $E(G^{-1})$ and the correspondence
structure on $E^{\sigma},$ as described in Proposition~\ref{corres},
fits the correspondence structure of $E(G^{-1})$. Consequently, we
can identify the two and write \[
E^{\sigma}=E(G^{-1}).\]
 (See Example 4.3 in \cite{MSHardy} for a description of the structure
of the dual correspondence for more general representations $\sigma$). It will also be convenient to write $\eta$ matricially with respect
to the orthonormal bases $\{\delta_{v}\mid v\in G^0\}$ of $H_{0}$
and $\{\delta_{e}\otimes e_{s(e)}\}_{e\in G^1}$ of $E\otimes H_{0}$
as \begin{equation}
(\eta)_{e,r(e)}=\eta(e^{-1}).\label{matrix1}\end{equation}

\end{example}

\subsection{Representations of the Hardy algebras}

We now turn to study the ultra-weakly continuous, completely contractive
representations of the Hardy algebra $H^{\infty}(E)$. Given such
a representation $\rho$, its restriction to the tensor algebra $\mathcal{T}_{+}(E)$
is a completely contractive representation of this algebra. Thus it
is of the form $T\times\sigma$ for some $\tilde{T}^{*}\in E^{\sigma}$
(where $\sigma=\rho\circ\varphi_{\infty}$). Therefore, the problem
we face is to decide when the integrated form, $T\times\sigma$, of
a completely contractive covariant representation $(T,\sigma)$ extends
from $\mathcal{T}_{+}(E)$ to $H^{\infty}(E)$. This problem arises
already in the simplest situation, vis. when $M=\mathbb{C}=E$. In
this setting, $T$ is given by a single contraction operator on a
Hilbert space, $\mathcal{T}_{+}(E)$ {}``is\textquotedblright{}\ the
disc algebra $A(\mathbb{D})$ and $H^{\infty}(E)$ {}``is\textquotedblright{}\ the
space $H^{\infty}(\mathbb{D})$ of bounded analytic functions on the
disc. In this case it is known that the representation $T\times\sigma$
extends from the disc algebra to $H^{\infty}(\mathbb{D})$ precisely
when there is no singular part to the spectral measure of the minimal
unitary dilation of $T$. In our general context, one may be able
to identify an analogue of a unitary dilation but it is rarely unique
(\cite{MS98}). Also, it doesn't seem to have any analogue for a spectral
measure. Thus we will need to use different tools.

One class of representations of the tensor algebra that extend to
ultra-weakly continuous representations of $H^{\infty}(E)$ we have
already met. These are the induced representations. In the notation
of (\ref{induced}), the ultra-weakly continuous extension of $\pi_{0}^{\mathcal{F}(E)}|\mathcal{T}_{+}(E)$
is $\pi_{0}^{\mathcal{F}(E)}|H^{\infty}(E)$.

It was proved in \cite[Theorem 2.13]{MSHardy} that, if $\Vert\tilde{T}\Vert<1$,
then the minimal isometric dilation $(V,\tau)$ of $(T,\sigma)$ (as
in Theorem~\ref{isomdil}) is an induced representation. Thus $V\times\tau$
extends to an ultra-weakly continuous representation of $H^{\infty}(E)$.
Since $T\times\sigma$ is a compression of $V\times\tau$, we have
the following.

\begin{lemma} \label{contraction} \cite[Corollary 2.14]{MSHardy}
If $\Vert\tilde{T}\Vert<1$ then $T\times\sigma$ extends to a ultraweakly
continuous representation of $H^{\infty}(E)$. \end{lemma}

If $T\times\sigma$ is a representation of the tensor algebra on the
space $H$ that extends to a ultra-weakly continuous representation
of $H^{\infty}(E)$ then, for every $x\in H$, the linear functional
$f=\omega_{x}\circ(T\times\sigma)$ extends to a ultra-weakly continuous
functional on $H^{\infty}(E)$. Given an arbitrary representation
$T\times\sigma$, one can still consider the set of all vectors $x\in H$
with this property. For the case where $M=\mathbb{C}$, this was done
in \cite{DLP} and the following definition is a direct extension
of their definition.

\begin{definition} Given a c.c. covariant representation $(T,\sigma)$ on $H$,
we say that $x\in H$ is \textbf{absolutely continuous} if the functional
$\omega_{x}\circ(T\times\sigma)$, on $\mathcal{T}_{+}(E)$, extends
to a ultraweakly continuous functional on $H^{\infty}(E)$ and we
write \textbf{$\mathcal{V}_{ac}(T,\sigma)$} for the set of all the
absolutely continuous vectors for $(T,\sigma)$. \end{definition}

It turns out that the set of absolutely continuous vectors can be
studied by considering the ranges of certain intertwiners.

\begin{definition} Let $(S_{0},\sigma_{0})$ be the universal induced
covariant representation (see Definition~\ref{Definition:Universal}).
For a given $\eta\in\overline{\mathbb{D}(E^{\sigma})}$ (corresponding
to the representation $(T,\sigma)$ on $H$) write $\mathcal{I}(S_{0},\eta^{*})$
(or $\mathcal{I}(S_{0},\tilde{T})$) for the space of intertwiners:
$\mathcal{I}(S_{0},\eta^{*})=$
\[
\{C:H_{0}\rightarrow H:CS_{0}(\xi)=T(\xi)C,\; C\sigma_{0}(a)=\sigma(a)C,\xi\in
E,\; a\in M\}.\]
 \end{definition}

The role of these intertwiners for studying the ultraweakly continuous
representations of $H^{\infty}(E)$ is seen in \cite[Corollaries 5.4 and 5.5]{D69}
(in the case $M=\mathbb{C}$) and in \cite[Lemma 7.12]{MSHardy} (in
the general case). An immediate corollary of the latter lemma is that,
when $(T,\sigma)$ is an isometric representation, for every $C\in\mathcal{I}(S_{0},\eta^{*})$,
the range $Ran(C)$ of $C$ is contained in $\mathcal{V}_{ac}(T,\sigma)$.
(Here $\eta=\tilde{T}^{*}$.) Generalizing some arguments of \cite[Theorem 1.6]{DLP},
we prove in \cite[Proposition 3.5]{MSAbsCont} that the converse also
holds. Thus we have the following.

\begin{theorem} If $(T,\sigma)$ is an isometric covariant representation then
\[
\mathcal{V}_{ac}(T,\sigma)=\bigcup\{Ran(C):C\in\mathcal{I}(S_{0},\tilde{T})\}.\]
 It follows that $\mathcal{V}_{ac}(T,\sigma)$ is a closed, $\sigma(M)$-invariant
subspace. \end{theorem} 


In order to analyze the set of all absolutely continuous vectors for
a general completely contractive (not necessarily isometric) representation,
we use the fact that every such representation $(T,\sigma)$ on $H$
has a minimal isometric dilation to an isometric representation $(V,\tau)$
on a larger space $K$ (Theorem~\ref{isomdil}). This dilation is
constructed explicitly in \cite{MS98} and it is evident from the
construction that the restriction of $(V,\tau)$ to $K\ominus H$
is an induced representation. It then follows that $K\ominus H\subseteq\mathcal{V}_{ac}(V,\tau)$.
Another tool used in the proof of the following theorem, which shows
the close relationship between $\mathcal{V}_{ac}(V,\tau)$ and $\mathcal{V}_{ac}(T,\sigma)$,
is the commutant lifting theorem (Theorem~\ref{CLT}). This theorem
is applied to show that every operator $C$ in $\mathcal{I}(S_{0},\tilde{T})$
can be written as $P_{H}X$ for some $X\in\mathcal{I}(S_{0},\tilde{V})$
(where $P_{H}$ is the projection onto $H$).

\begin{theorem}\label{characterization} Let $(T,\sigma)$ be a completely
contractive covariant representation of $(E,M)$ on the Hilbert space $H$,
let $(V,\rho)$ be the minimal isometric dilation of $(T,\sigma)$
acting on a Hilbert space $K$ containing $H$, and let $P$ denote
the projection of $K$ onto $H$. Then $K\ominus H$ is contained
in $\mathcal{V}_{ac}(V,\rho)$ and the following sets are equal.
\begin{enumerate}
\item [(1)] $\mathcal{V}_{ac}(T,\sigma)$.
\item [(2)]$H\cap\mathcal{V}_{ac}(V,\rho)$.
\item [(3)] $P\mathcal{V}_{ac}(V,\rho)$.
\item [(4)]$\bigcup\{Ran(C)\mid C\in\mathcal{I}(S_{0},\tilde{T})\}$.

\end{enumerate}
In particular, $\mathcal{V}_{ac}(T,\sigma)=H$ if and only if $\mathcal{V}_{ac}(V,\rho)=K$.

\end{theorem}

\begin{definition} Given a c.c. covariant representation $(T,\sigma)$, one
defines \textbf{the completely positive map associated with it}, $\Phi_{T}:\sigma(M)'\rightarrow\sigma(M)'$
by \[
\Phi_{T}(b)=\eta^{*}(I_{E}\otimes b)\eta=\tilde{T}(I_{E}\otimes b)\tilde{T}^{*}\]
 where $\eta=\tilde{T}^{*}$. \end{definition}

\begin{example} If $M=\mathbb{C}$ and $E=\mathbb{C}^{d}$, every
completely contractive representation is given by a row contraction
$\tilde{T}=(T_{1},\ldots,T_{d})$. In this case \[
\Phi_{T}(b)=\sum T_{i}bT_{i}^{*}.\]
 \end{example}

Note that the map $\Phi_{T}$ is a completely positive, contractive,
normal map on the von Neumann algebra $\sigma(M)'$. In fact, \textbf{every}
completely positive, contractive, normal map on a von Neumann algebra
$N$ is of this form. (See \cite[Corollary 2.23]{MSQMP} for details).

Often properties of the representation $(T,\sigma)$ can be expressed
in terms of the associated map $\Phi_{T}$. For example, the map is
multiplicative (that is, a $^{*}$-endomorphism) if and only if the
representation is isometric. (This is a rough statement. For the precise
one see \cite[Proposition 2.21]{MSQMP}). Another example is the curvature
associated with a representation which was shown in \cite{MSCurv}
to be an artifact of the associated map.

Here, too, we find that there is a close relationship between the
intertwiners in $\mathcal{I}(S_{0},\eta^{*})$ and pure superharmonic
elements (to be defined below) of the associated map.


\begin{lemma} If $C\in\mathcal{I}(S_{0},\eta^{*})$, then $Q=CC^{*}$
lies in $\sigma(M)'$ and satisfies
\begin{itemize}
\item [(i)] $Q\geq0$ and $\Phi_{T}(Q)\leq Q$, and
\item [(ii)] $\Phi_{T}^{n}(Q)\rightarrow0$ ultra weakly.
\end{itemize}
\end{lemma} The proof of (i) follows from the fact that $CS_{0}(\xi)=T(\xi)C$
for $\xi\in E$ and, thus, $C\tilde{S_{0}}=\tilde{T}(I_{E}\otimes C)$
and $\Phi_{T}(CC^{*})=\tilde{T}(I_{E}\otimes CC^{*})\tilde{T}^{*}=C\tilde{S_{0}}\tilde{S_{0}}^{*}C^{*}\leq CC^{*}$.
The proof of (ii) follows similarly from the fact that $S_{0}$ is
pure (that is, $\tilde{S_{0}}_{n}\tilde{S_{0}}_{n}^{*}\rightarrow0$
strongly as $n\to\infty$).

\begin{definition} An element $Q\in\sigma(M)'$ satisfying (i) of
the lemma will be said to be \emph{ superharmonic} for $\Phi_{T}$.
If it also satisfied (ii), it will be said to be \emph{pure superharmonic}.
\end{definition}

Thus, for every $C\in\mathcal{I}(S_{0},\eta^{*})$, $Q=CC^{*}$ is
pure superharmonic. The converse also holds. Given an element $Q\in\sigma(M)'$
that is pure superharmonic for $\Phi_{T}$, define $r\in\sigma(M)'$
to be the positive square root of $Q-\Phi_{T}(Q)$. Then $\sum_{n\geq0}\Phi_{T}^{n}(r^{2})=Q$
(in the strong operator topology).

Let $\mathcal{R}$ be the closure of the range of $r$ and $\sigma_{\mathcal{R}}:=\sigma|\mathcal{R}$.
Since $\pi$ is a normal, faithful representation of $M$ with infinite
multiplicity and $\sigma_{\mathcal{R}}$ is a normal representation
of $M$, there exists an isometry $v:\mathcal{R}\rightarrow K_{0}$
that intertwines $\sigma_{\mathcal{R}}$ and $\pi$. Then set \[
C^{*}=(I_{\mathcal{F}(E)}\otimes v)\sum_{n\geq0}(I_{E^{\otimes n}}\otimes r)\tilde{T}_{n}^{*}\]
 where $\tilde{T}_{n}=\tilde{T}(I_{E}\otimes\tilde{T})\cdots(I_{E^{\otimes(n-1)}}\otimes\tilde{T}):E^{\otimes n}\otimes_{\sigma}H\rightarrow H$.

It is then straightforward to show that $C\in\mathcal{I}(S_{0},\tilde{T})$
and $CC^{*}=Q$. Some of the arguments above can be found, for $M=\mathbb{C}$,
in \cite{D69} and in \cite{gP03}. We conclude the following. \begin{corollary}\label{AC_and_intertwiners}
\[
\mathcal{V}_{ac}(T,\sigma)=\bigcup\{Ran(C):C\in\mathcal{I}(S_{0},\tilde{T})\}=\]
 \[
=\bigvee\{Ran(Q):Q\mbox{ is a pure superharmonic operator for }\Phi_{T}\}.\]
 \end{corollary}

We can now state the following result which gives a complete description
of the representations of $H^{\infty}(E)$ (see \cite[Theorem 4.1]{MSAbsCont}).

\begin{theorem}\label{AbsContConditions} Let $T\times\sigma$ be
a c.c. representation of $\mathcal{T}_{+}(E)$ on $H$ and write $\eta=\tilde{T}^{*}$
for the element of $\overline{\mathbb{D}(E^{\sigma})}$ associated
with it. Then the following are equivalent.
\begin{itemize}
\item [(1)] The representation $T\times\sigma$ extends to a completely
contractive ultra weakly continuous representation of $H^{\infty}(E)$.
\item [(2)] $\mathcal{V}_{ac}(T,\sigma)=H$
\item [(3)]$H=\bigvee\{Ran(C):C\in\mathcal{I}(S_{0},\eta^{*})\}$.
\item [(4)] $H=\bigvee\{Ran(Q):Q\;\;\mbox{is pure superharmonic for}\;\Phi_{T}\}$
\end{itemize}
\end{theorem}

Note that the equivalence of (1) and (2) of the theorem means that
this {}``extension\textquotedbl{} problem can be studied {}``locally\textquotedbl{}.

Theorem~\ref{AbsContConditions} describes the representations $\rho$
of $H^{\infty}(E)$ that satisfy $\rho\circ\varphi_{\infty}=\sigma$.
The set of the points $\eta\in\overline{\mathbb{D}(E^{\sigma})}$
that correspond to these representations will be denoted $AC(E^{\sigma})$.
These sets (for all $\sigma$'s) parameterize the representations
of $H^{\infty}(E)$. We have \[
\mathbb{D}(E^{\sigma})\subseteq AC(E^{\sigma})\subseteq\overline{\mathbb{D}(E^{\sigma})}.\]

As corollaries of the analysis above, we get the following.

\begin{theorem} \cite[Theorem 5.6]{MSAbsCont} If $\sigma(M)'$ is
finite dimensional, then $\rho:=T\times\sigma$ extends to an ultra-weakly
continuous representation of $H^{\infty}(E)$ on $H$ if and only
if $(T,\sigma)$ is completely non coisometric; that is, there is
no $\rho(\mathcal{T}_{+}(E))^{*}$-invariant subspace of $H$ on which
$\tilde{T}^{*}$ is an isometry. \end{theorem}

\begin{theorem} \cite[Theorem 5.3]{MSAbsCont} If $\sigma(M)'$ has
a non zero normal periodic state $\omega$ for $\Phi_{T}$ (that is,
$\omega\circ\Phi_{T}^{k}=\omega$ for some $k\geq1$) then $T\times\sigma$
does not extend to an ultra-weakly continuous representation of $H^{\infty}(E)$.
\end{theorem}

When $(T,\sigma)$ is an isometric representation, the space $\mathcal{V}_{ac}(T,\sigma)$
contains all the wandering vectors of $T\times\sigma$ where $h\in H$
is said to be wandering if, for every $n\neq m$, the subspaces $\tilde{T}_{n}(E^{\otimes n}\otimes[\sigma(M)h])$
and $\tilde{T}_{m}(E^{\otimes m}\otimes[\sigma(M)h])$ are orthogonal.
This, and other results concerning wandering vectors, was proved in
\cite{MSAbsCont} generalizing similar results in \cite{DLP} (who
proved it for the case $M=\mathbb{C}$).

\section{The functions defined by elements of $H^{\infty}(E)$}

As we stated in the introduction, we view the elements of $H^{\infty}(E)$
as functions defined on the space of the representations of $H^{\infty}(E)$.

As seen in the previous section, the space of all representations
can be parameterized by $\cup AC(E^{\sigma})^{*}$ (where the union
runs over all normal representations $\sigma$ of $M$). In the discussion
below, we shall fix $\sigma$.

Now, a given $X\in H^{\infty}(E)$ will be viewed as a function $\widehat{X}$,
defined on $AC(E^{\sigma})^{*}$ by the equation \begin{equation}
\widehat{X}(\eta^{*})=(\eta^{*}\times\sigma)(X)\;,\quad\eta\in AC(E^{\sigma}).\label{pointevaluation}\end{equation}

Often it will be more convenient to restrict the function $\widehat{X}$
to the open unit ball $\mathbb{D}(E^{\sigma})^{*}$. This restriction
will also be denoted $\widehat{X}$.

The primary objective in this section is to understand the range of
the transform \[
X\mapsto\widehat{X}\]
 from $H^{\infty}(E)$ to the set of all $B(H)$-valued functions
on $AC(E^{\sigma})^{*}$ or on $\mathbb{D}(E^{\sigma})^{*}$.

Before we do this we note that this map depends on $\sigma$ and that,
for a given $\sigma$, it may have a non zero kernel; that is, there
may be some $X\in H^{\infty}(E)$ such that $\widehat{X}=0$. However,
it was shown in \cite[Lemma 5.7]{MSSchur} that we can always choose
an appropriate $\sigma$ so that this transform is injective. We shall
have more to say about the kernel of the transform in Section 5.


\begin{example} \label{ex0.12}Suppose $M=E=\mathbb{C}$ and $\sigma$
the representation of $\mathbb{C}$ on some Hilbert space $H$. Then
it is easy to check that $E^{\sigma}$ is isomorphic to $B(H)$. Fix
an $X\in H^{\infty}(E)$. As we mentioned above, this Hardy algebra
is the classical $H^{\infty}(\mathbb{D})$ and we can identify $X$
with a function $f\in H^{\infty}(\mathbb{T})$. Given $S\in\mathbb{D}(E^{\sigma})=B(H)$,
it is not hard to check that $\widehat{X}(S^{\ast})$, as defined
above, is the operator $f(S^{\ast})$ defined through the usual $H^{\infty}$-functional
calculus. \end{example}

\begin{example} \label{freesgp}In \cite{DP98} Davidson and Pitts
associate to every element of the algebra $\mathcal{L}_{n}=H^{\infty}(\mathbb{C}^{n})$
a function on the open unit ball of $\mathbb{C}^{n}$. This is a special
case of our analysis when $M=\mathbb{C}$, $E=\mathbb{C}^{n}$ and
$\sigma$ is a one dimensional representation of $\mathbb{C}$. In
this case $\sigma(M)^{\prime}=\mathbb{C}$ and $E^{\sigma}=\mathbb{C}^{n}$.
Note, however, that our definition allows us to take $\sigma$ to
be the representation of $\mathbb{C}$ on an arbitrary Hilbert space
$H$. If we do so, then $E^{\sigma}$ is isomorphic to $B(H)^{(n)}$,
the nth column space over $B(H)$, and elements of $\mathcal{L}_{n}$
define functions on the open unit ball of this space viewed as a correspondence
over $B(H)$ with values in $B(H)$. 
\end{example}

Note that, if $\xi_{n}\in E^{\otimes n}$ and $X=T_{\xi_{n}}$, we
have, for $\eta\in AC(E^{\sigma})$ and $h\in H$, \[
\widehat{T_{\xi_{n}}}(\eta^{*})h=\eta_{n}^{*}(\xi_{n}\otimes h)\]
 where, recall, $\eta_{n}=(I_{E^{\otimes(n-1)}}\otimes\eta)\cdots(I_{E}\otimes\eta)\eta:H\rightarrow E^{\otimes n}\otimes H$.

\begin{example}\label{quiverevaluation} Let $G$ and $E=E(G)$ be
as in Example~\ref{quivercorrespondence}. Let $\{e_{v}:v\in G^{0}\}$
be an orthonormal basis for a Hilbert space $H$ and let $\sigma$
be the diagonal representation of $\ell^{\infty}(G^{0})$ on $H$
(as at the end of Example~\ref{quiverrepresentations}). It follows
from the computations in that example that, for the generators of
$H^{\infty}(G)$, we get \begin{equation}
\widehat{P_{v}}(\eta^{\ast})=\theta_{v,v}\;,\; v\in V\label{eval1}\end{equation}
 and \begin{equation}
\widehat{S_{e}}(\eta^{\ast})=\overline{\eta(e^{-1})}\theta_{r(e),s(e)}\;,\; e\in\mathcal{Q}\label{eval2}\end{equation}
 where $\theta_{v,w}$ is the partial isometry operator on $H$ that
maps $e_{w}$ to $e_{v}$ and vanishes on $(e_{w})^{\perp}$. For
a general $X\in H^{\infty}(G)$, $\widehat{X}(\eta^{\ast})$ is obtained
by using the linearity, multiplicativity and $w^{\ast}$-continuity
of the map $X\mapsto\widehat{X}(\eta^{\ast})$.

\end{example}

In order to understand what functions can be obtained as $\widehat{X}$
for some $X\in H^{\infty}(E)$, we first consider the following simple
examples.

\begin{example} Let $M=E=\mathbb{C}$ and $\sigma$ be the one dimensional
representation of $\mathbb{C}$. The Hardy algebra is the classical
$H^{\infty}(\mathbb{D})$ and, for $X\in H^{\infty}(\mathbb{D})$,
the function $\widehat{X}$ is just $X$. Thus the functions we get
are the functions on $\mathbb{D}$ (which is $\mathbb{D}(E^{\sigma})$
in this case) that are holomorphic and bounded. Recall that the Schur
class $\mathcal{S}$ is the set of all such functions $S$ with $|S(z)|\leq1$
(for $z\in\mathbb{D}$). \end{example} \begin{example} Let $M=E=B(H)$
and $\sigma$ be the identity representation of $B(H)$ on $H$. The
Fock space $\mathcal{F}(E)$ is the direct sum of infinitely many
copies of $B(H)$ and $\mathcal{F}(E)\otimes_{\sigma}H$ is (isomorphic
to) $\ell^{2}\otimes H$. The operator $T_{I}\otimes I_{H}$ is $S\otimes I_{H}$
where $S$ is the unilateral shift on $\ell^{2}$ and the operator
$\varphi_{\infty}(A)\otimes I_{H}$ (for $A\in M=B(H)$) is $I_{\ell^{2}}\otimes A$.
Since the induced representation $\sigma^{\mathcal{F}(E)}|H^{\infty}(E)$
is completely isometric and a homeomorphism with respect to the
ultra-weak topologies, we can identify $H^{\infty}(E)$ with $H^{\infty}(\mathbb{D})\otimes B(H)$.
The $\sigma$-dual $E^{\sigma}$ in this case is $\mathbb{C}$ and
$\mathbb{D}(E^{\sigma})=\mathbb{D}$. The transform $X\mapsto\widehat{X}$
is just the expression of an element in $H^{\infty}(\mathbb{D})\otimes B(H)$
as a bounded $B(H)$-valued holomorphic function on $\mathbb{D}$.
The set of all such functions $S$ that satisfy $\norm{S(z)}\leq1$
(for all $z\in\mathbb{D}$) is known as the operator valued Schur
class $\mathcal{S}(H)$. \end{example}

Our objective in this section is to show that, in the general case,
the functions we get as $\widehat{X}$ (for $X\in H^{\infty}(E)$
with $\norm{X}\leq1$ ) should be viewed as generalized Schur class
functions.

For this, recall first, that the functions in the operator valued
Schur class have several characterizations. The following is well
known (see \cite{BBFH} for a more detailed exposition of the operator-valued
Schur class and some of its generalizations).

\begin{theorem}\label{schurclass} For an $B(H)$-valued function
$S$ on $\mathbb{D}$ the following conditions are equivalent.
\begin{itemize}
\item [(1)] $S\in\mathcal{S}(H)$; that is, $S$ is a $B(H)$-valued holomorphic
function on $\mathbb{D}$ with $\norm{S(z)}\leq1$ for all $z\in\mathbb{D}$.
\item [(2)] There is a Hilbert space $\mathcal{E}$ and a coisometric
operator ( called colligation) \[
U=\left(\begin{array}{cc}
A & B\\
C & D\end{array}\right):\left(\begin{array}{c}
\mathcal{E}\\
H\end{array}\right)\rightarrow\left(\begin{array}{c}
\mathcal{E}\\
H\end{array}\right)\]
 so that $S$ can be realized as a linear fractional function \begin{equation}
S(z)=D+zC(I_{\mathcal{E}}-zA)^{-1}B.\label{transfer}\end{equation}

\item [(3)] The function $K_{S}:\mathbb{D}\times\mathbb{D}\rightarrow B(H)$
given by

\[
K_{S}(z,w)=\frac{I-S(z)S(w)^{*}}{1-z\overline{w}}\]
 is a positive kernel on $\mathbb{D}\times\mathbb{D}$ (with values
in $B(H)$).

\end{itemize}
\end{theorem}

From the point of view of systems theory, a function $S$ realized
as in (\ref{transfer}) is the transfer function of a certain linear
system (defined using $A,B,C$ and $D$).

For (3) above, recall that a function $K:\Omega\times\Omega\rightarrow B(H)$
is said to be a positive kernel if, for every $k\geq1$ and every
choice of $\omega_{1},\ldots,\omega_{k}$ in $\Omega$, the matrix
$(K(\omega_{i},\omega_{j}))$ is positive (as an element of $M_{k}(B(H))$.

In order to discuss generalized Schur class operator functions we
need to define a completely positive definite kernel. The definition
below can be found in \cite[Definition 3.2.2]{BBLS}. In fact, in
that paper the authors show that this definition is equivalent to
several other definitions and they prove an extension of Kolmogorov's
representation theorem for these kernels.

\begin{definition} Let $A,C$ be $C^{*}$-algebras and write $B(A,C)$
for the bounded (linear) maps from $A$ to $C$. Let $\Omega$ be
an arbitrary set. A function $K:\Omega\times\Omega\rightarrow B(A,C)$
is said to be a completely positive definite kernel if, for every
$k\geq1$ and every choice of $\omega_{1},\ldots,\omega_{k}$ in $\Omega$,
the map $\Psi_{K}:M_{k}(A)\rightarrow M_{k}(C)$ defined by $\Psi_{K}((a_{i,j}))=(K(\omega_{i},\omega_{j})(a_{i,j}))$
is completely positive. \end{definition}

The following two theorems (Theorem~\ref{Realization} and Theorem~\ref{kernel}),
when compared with Theorem~\ref{schurclass}, show that we can indeed
view the functions $\{\widehat{X}:X\in H^{\infty}(E)\;,\;\;\norm{X}\leq1\}$
as Schur class operator functions. The proofs can be found in \cite{MSSchur}.

\begin{theorem}\label{Realization} Let $E$ be a $W^{*}$-correspondence
over $M$, $\sigma$ a faithful normal representation of $M$ on $H$
and $Z:\mathbb{D}(E^{\sigma})^{*}\rightarrow B(H)$. Then $Z=\hat{X}$
for some $X\in H^{\infty}(E)$ with $\norm{X}\leq1$ if and only if
there is a Hilbert space $\mathcal{E}$, a normal representation $\tau$
of $\sigma(M)'$ on $\mathcal{E}$ and a coisometric operator matrix
\[
U=\left(\begin{array}{cc}
A & B\\
C & D\end{array}\right):\left(\begin{array}{c}
\mathcal{E}\\
H\end{array}\right)\rightarrow\left(\begin{array}{c}
E^{\sigma}\otimes_{\tau}\mathcal{E}\\
H\end{array}\right)\]
 (with $A,B,C,D$ that are $\sigma(M)'$-module maps ) so that $Z$
can be realized as \[
Z(\eta^{*})=D+C(I_{\mathcal{E}}-L_{\eta}^{*}A)^{-1}L_{\eta}^{*}B.\]
 \end{theorem} Here $L_{\eta}:\mathcal{E}\rightarrow E^{\sigma}\otimes_{\tau}\mathcal{E}$
is defined by $L_{\eta}h=\eta\otimes h$.

\begin{theorem}\label{kernel} Let $E$ be a $W^{*}$-correspondence
over $M$, $\sigma$ a faithful normal representation of $M$ on $H$
and $Z:\mathbb{D}(E^{\sigma})^{*}\rightarrow B(H)$. Then $Z=\hat{X}$
for some $X\in H^{\infty}(E)$ with $\norm{X}\leq1$ if and only if
the kernel $K_{Z}:\mathbb{D}(E^{\sigma})^{*}\times\mathbb{D}(E^{\sigma})^{*}\rightarrow B(\sigma(M)',B(H))$
is completely positive definite where \[
K_{Z}(\eta^{*},\zeta^{*})=(id-Ad(Z(\eta^{*}),Z(\zeta^{*})))\circ(id-\theta_{\eta,\zeta})^{-1}.\]
 \end{theorem} Here $Ad(Z(\eta^{*}),Z(\zeta^{*}))(a)=Z(\eta^{*})aZ(\zeta^{*})^{*}$
and $\theta_{\eta,\zeta}(a)=\langle\eta,a\zeta\rangle$ for $a\in\sigma(M)'$.

In \cite{MSSchur} we used the condition appearing in Theorem~\ref{kernel}
to define a Schur class operator function.

\begin{definition} \label{schur}Let $\Omega$ be a subset of $\mathbb{D}(E^{\sigma})$
and let $\Omega^{\ast}=\{\omega^{\ast}\mid\omega\in\Omega\}$. A function
$Z:{\Omega}^{\ast}\rightarrow B(H)$ will be called a \emph{Schur
class operator function} (with values in $B(H)$) if the kernel $K_{Z}:\Omega^{*}\times\Omega^{*}\rightarrow B(\sigma(M)',B(H))$,
defined by \[
K_{Z}(\eta^{*},\zeta^{*})=(id-Ad(Z(\eta^{*}),Z(\zeta^{*})))\circ(id-\theta_{\eta,\zeta})^{-1}\;,\;\;\;\eta,\zeta\in\Omega\]
 is completely positive definite. 
\end{definition}

Thus, we see that the functions of the form $\widehat{X}$ are precisely
the Schur class operator functions on the open unit ball of $E^{\sigma}$.

One direction of Theorem~\ref{kernel} follows from a Nevanlinna-Pick
type interpolation theorem that we proved in \cite[Theorem 5.3]{MSHardy}.
Before we state it, recall the classical Nevanlinna-Pick theorem.

\begin{theorem} Given $z_{1},\ldots,z_{m}$ in $\mathbb{D}$ and
$w_{1},\ldots,w_{m}$ in $\mathbb{C}$, one can find a function $f\in H^{\infty}(\mathbb{D})$
with $\norm{f}\leq1$ and $f(z_{i})=w_{i}$ for all $i$ if and only
if the $m\times m$ matrix \[
\left(\frac{1-w_{i}\overline{w_{j}}}{1-z_{i}\overline{z_{j}}}\right)\]
 is positive. \end{theorem}

Our generalization is Theorem~\ref{NP}. This result captures numerous
theorems in the literature that go under the name of generalized Nevanlinna-Pick
theorems. In particular, of course, it gives the classical Nevanlinna-Pick
theorem (when $M=E=\mathbb{C}$). In the setting when $M=\mathbb{C}$
and $E=\mathbb{C}^{n}$, it gives versions due to Popescu in \cite{gP98},
Arias and Popescu in \cite{AP00} and Davidson and Pitts in \cite{DP98b}.

\begin{theorem} \label{NP}(\cite[Theorem 5.3]{MSHardy}) Let $E$
be a $W^{\ast}$-correspondence over a von Neumann algebra $M$ and
let $\sigma:M\rightarrow B(H)$ be a faithful normal representation
of $M$ on a Hilbert space $H$. Fix $k$ points $\eta_{1},\ldots\eta_{k}\ $in
the disk $\mathbb{D}(E^{\sigma})$ and choose $2k$ operators $B_{1},\ldots B_{k},C_{1},\ldots C_{k}$
in $B(H)$. Then there exists an $X$ in $H^{\infty}(E)$ such that
$\Vert X\Vert\leq1$ and\[
B_{i}\widehat{X}(\eta_{i}^{\ast})=C_{i}\]
 for $i=1,2,\ldots,k,$ if and only if the map from $M_{k}(\sigma(M)^{\prime})$
to $M_{k}(B(H))$ defined by the $k\times k$ matrix \begin{equation}
\left((Ad(B_{i},B_{j})-Ad(C_{i},C_{j}))\circ(id-\theta_{\eta_{i},\eta_{j}})^{-1}\right)\label{PickMatrix}\end{equation}
 is completely positive. \end{theorem}

The proof was inspired by Sarason's approach that is based on the
commutant lifting theorem \cite{dS67} and the extensions of it to
the multi-analytic setting of \cite{gP98}. The organization of our
proof follows the presentation of \cite{RR97} and one of the key
ingredients is our commutant lifting theorem (Theorem~\ref{CLT}).
Another important ingredient in the proof is the following theorem
that identifies the commutant of an induced representation and points
out another role that the $\sigma$-dual $E^{\sigma}$ plays in studying
the Hardy algebra $H^{\infty}(E)$.

In order to state the theorem, recall that if $E$ is a $W^{*}$-correspondence
over the von Neumann algebra $M$ and $\sigma$ is a faithful normal
representation of $M$ on $H$, we can represent $H^{\infty}(E)$
on $\mathcal{F}(E)\otimes_{\sigma}H$ using the induced representation.
Write $\rho$ for this representation so that $\rho(X)=\sigma^{\mathcal{F}(E)}(X)=X\otimes I_{H}$
for $X\in H^{\infty}(E)$. Note that this representation is completely
isometric isomorphism and a homeomorphism with respect to the ultra-weak
topologies. Similarly, we have an induced representation $\rho'$
of $H^{\infty}(E^{\sigma})$ on $\mathcal{F}(E^{\sigma})\otimes_{\iota}H$
where $\iota$ is the identity representation of $\sigma(M)'$ on
$H$.

\begin{theorem}\label{commutant} (\cite[Theorem 3.9]{MSHardy})
Let $E$ be a $W^{*}$-correspondence over the von Neumann algebra
$M$, let $\sigma$ be a faithful normal representation of $M$ on
$H$ and let $\rho$ and $\rho'$ be as above. Then the commutant
of $\rho(H^{\infty}(E))$ is unitarily isomorphic to $\rho'(H^{\infty}(E^{\sigma}))$.

Consequently (using duality arguments), $(\rho(H^{\infty}(E)))''=\rho(H^{\infty}(E))$. \end{theorem}

Using the Nevanlinna-Pick theorem, we get another interesting result
that fits with the {}``noncommutative function theory\textquotedbl{}
point of view. This is the generalization of Schwartz's lemma (see
\cite[Theorem 5.6]{MSHardy}). It asserts, among several things, that
for $X\in H^{\infty}(E)$, if $\widehat{X}$ vanishes at the origin and if $\left\Vert X\right\Vert \leq1$,
then \[
\widehat{X}(\eta^{\ast})\widehat{X}(\eta^{\ast})^{\ast}\leq\langle\eta,\eta\rangle\text{,}\quad\eta\in\mathbb{D}(E^{\sigma})\]
 where, recall, $\langle\cdot,\cdot\rangle$ is the $\sigma(M)^{\prime}$-valued
inner product on $E^{\sigma}$ defined above. (See Proposition~\ref{corres}).

The generalized Nevanlinna-Pick theorem (Theorem~\ref{NP}), as stated
above, interpolates the values of the function at points in the open
unit ball $\mathbb{D}(E^{\sigma})^{*}$. But we now know that functions
of the form $\widehat{X}$ are defined on the, possibly larger, set
$AC(E^{\sigma})^{*}$. We wish to present a similar interpolation
theorem where the points $\eta_{1},\ldots,\eta_{k}$ are from $AC(E^{\sigma})$.

The first problem that one encounters in trying to do this is that
the maps $(id-\theta_{\eta_{i},\eta_{j}})$, appearing in the theorem,
are not necessarily invertible (as we may have $\norm{\eta_{i}}=1$).
The theorem, therefore, will have to be stated differently. In order
to deal with points on the boundary we use the following simple observation.

\textbf{Simple observation :} If $\Phi,\Psi$ are positive maps such
that $id-\Phi$ is invertible then the map $(id-\Psi)\circ(id-\Phi)^{-1}$
is positive if and only if : \[
\{a\geq0:\Phi(a)\leq a\}\subseteq\{a\geq0:\Psi(a)\leq a\}.\]

The last statement makes sense even if $id-\Phi$ is not invertible.
It is related to the Lyapunov preorder studied in matrix theory. This
was pointed out to us by Nir Cohen.

Using this observation and the characterizations of the points in
$AC(E^{\sigma})$ (Theorem~\ref{characterization}) we get the following.

\begin{theorem} Let $E$ be a $W^{*}$-correspondence over the von
Neumann algebra $M$ and let $\sigma$ be a faithful normal representation
of $M$ on $H$. Given $\eta_{1},\eta_{2},\ldots,\eta_{k}\in AC(E^{\sigma})$
and $D_{1},D_{2},\ldots,D_{k}\in B(H)$, the following conditions
are equivalent.
\begin{enumerate}
\item [(1)] There is an element $X\in H^{\infty}(E)$ such that $\norm{X}\leq1$
and such that \[
\widehat{X}(\eta_{i}^{*})=D_{i}\;,\]
 $i=1,2,\ldots,k$.
\item [(2)] For each $m\geq1$, $i:\{1,\ldots,m\}\rightarrow\{1,\ldots,k\}$
and $C_{1},C_{2},\ldots,C_{m}$ with $C_{j}\in\mathcal{I}(S_{0},\eta_{i(j)}^{*})$,
we have \[
(D_{i(l)}C_{l}C_{j}^{*}D_{i(j)}^{*})_{l,j}\leq(C_{l}C_{j}^{*})_{l,j}.\]

\end{enumerate}
\end{theorem}

In the case where $M=E=\mathbb{C}$, the above theorem gives an answer the the question: Given $k$ contractions $T_1,T_2, \ldots,T_k$ in $B(H)$ that have $H^{\infty}$-functional calculus and $k$ operators $D_1,D_2,\ldots,D_k$ in $B(H)$, when can we find a function $h\in H^{\infty}(\mathbb{D})$ such that $D_i=h(T_i)$ for all $1\leq i \leq k$?

\section{The kernel of the transform $X\mapsto\widehat{X}$ and quotient algebras}

In the last section we discussed the map $X\mapsto\widehat{X}$ that
maps every element $X$ of $H^{\infty}(E)$ to a function $\widehat{X}$,
defined on $\mathbb{D}(E^{\sigma})^{*}$ or $AC(E^{\sigma})^{*}$,
that was seen to be a Schur class operator function. We noted there
that this transform depends on $\sigma$ and it may have a kernel.
This was observed already by Davidson and Pitts in \cite{DP98}. They
showed that, when $M=\mathbb{C}$, $E=\mathbb{C}^{n}$ and $\sigma$
is the one-dimensional representation of $\mathbb{C}$ (see Example~\ref{freesgp}),
the kernel of this map is the commutator ideal of $H^{\infty}(\mathbb{C}^{n})$.

In general, we write $K(\sigma)$ for this kernel. Thus \begin{equation}
K(\sigma)=\cap\{Ker(\eta^{*}\times\sigma):\;\;\eta\in\mathbb{D}(E^{\sigma})\}.\label{Ksig}\end{equation}

The following lemma was proved in \cite[Lemmas 3.7 and 4.17]{MSSchur}.

\begin{lemma} Let $\sigma$ be a normal, faithful, representation
of $M$ on a Hilbert space $H$ and let $K(\sigma)$ be defined by
(\ref{Ksig}). Then
\begin{enumerate}
\item [(1)] $K(\sigma)$ is an ultra-weakly closed ideal of $H^{\infty}(E)$.
\item [(2)] $K(\sigma)\subseteq\{X\in H^{\infty}(E):\;\Phi_{0}(X)=\Phi_{1}(X)=0\}$
(where $\Phi_{k}$ were defined in (\ref{FourierOperators})).
\item [(3)] $K(\sigma)$ is invariant under the action of the gauge group
and, thus, under the maps $\Phi_{k}$, $k\geq0$.
\item [(4)] If $\sigma$ is of infinite multiplicity, then $K(\sigma)=\{0\}$.
\end{enumerate}
\end{lemma}

It follows from statement (4) that we can always choose $\sigma$
such that this transform is injective.

For a general representation $\sigma$, we can view this transform
as a map on the quotient algebra $H^{\infty}(E)/K(\sigma)$ (into
the Schur class operator functions). Note, however, that, as shown by Arveson \cite{wA98}, this map is not isometric when the Schur class functions are viewed with the supremum norm. The following theorem identifies
quotient algebras of $H^{\infty}(E)$ with algebras that are obtained
by compressing $H^{\infty}(E)$ into a coinvariant subspace of $\mathcal{F}(E)$.
For the case $M=\mathbb{C}$ and $E=\mathbb{C}^{n}$, this was proved
by Davidson and Pitts in \cite{DP98b}. The general result was proved
independently by J. Meyer (\cite{jM}) and by M. Gurevich (\cite{mG}).
A norm-closed version was proved by A. Viselter (\cite{V}).

\begin{theorem} \label{quotient} Let $J\subseteq H^{\infty}(E)$
be an ultra-weakly closed two-sided ideal in $H^{\infty}(E)$. Write
$\mathcal{M}$ for the closed submodule $\overline{J\mathcal{F}(E)}$
and $P\in\mathcal{L}(\mathcal{F}(E))$ for the projection onto $\mathcal{M}^{\perp}$.
Then the map $X\in H^{\infty}(E)\mapsto PXP$ induces a complete isometric
isomorphism mapping the quotient algebra $H^{\infty}(E)/J$ onto $PH^{\infty}(E)P$.
\end{theorem}

When $J$ is invariant for the action of the gauge group on $H^{\infty}(E)$,
as is the case for the ideal $K(\sigma)$, the algebra $PH^{\infty}(E)P$
(for $P$ as in the theorem) is the Hardy algebra of a \emph{subproduct
system}. Subproduct systems were defined and studied in \cite{ShS}.
Theorem~\ref{quotient} shows that they can be useful in studying
the quotient algebras $H^{\infty}(E)/K(\sigma)$. We shall not discuss
it further here and the interested reader is referred to \cite{ShS}
and \cite{V}.

\section{Varying $\sigma$}

In most of the discussion above we fixed a normal representation $\sigma$
and, for $X\in H^{\infty}(E)$, considered the function $\widehat{X}$
defined on $\mathbb{D}(E^{\sigma})^{*}$ or on $AC(E^{\sigma})^{*}$.
Now we let $\sigma$ vary. We fix $M$ and $E$ and write $\Sigma$
for the set of all normal representations $\sigma$ of $M$ on some
Hilbert space $H_{\sigma}$. For every $X\in H^{\infty}(E)$ and every
$\sigma\in\Sigma$, we write $\widehat{X_{\sigma}}$ for the (Schur
class operator) function associated to $X$ on $AC(E^{\sigma})^{*}$.
We get a family of operator valued functions $\{\widehat{X_{\sigma}}:\sigma\in\Sigma\}$.

In this section we discuss the relationships among the functions in
this family. Our discussion was inspired by several sources. First,
there is the pioneering paper by Joe Taylor \cite{jT72}. This paper
seems to have generated very little interest until relatively recently.
But on close reading, it is clear that it was extraordinarily prescient.
It had a big impact on Dan Voiculescu's ``free analysis'' questions
\cite{Voi} and most recently it helped to shape the foundations of
noncommutative function theory being developed by Dimitry Kalyuzhny\u{i}-Verbovetzki\u{i}
and Victor Vinnikov and the applications of it to linear matrix inequalities
and real algebraic geometry in the work of Bill Helton, Igor Klepp
and Scott McCullough and their collaborators (see, e.g., \cite{HKM}).
And of course, it has been a direct source of inspiration for our
work.

Each function $\widehat{X_{\sigma}}$ is defined on its own domain,
$AC(E^{\sigma})^{*}$, and these domains vary with $\sigma$. One
can therefore view $\{\widehat{X_{\sigma}}\}_{\sigma\in\Sigma}$ as
a single function $\mathbf{\widehat{X}}$ from $\mathcal{AC}(E):=\coprod_{\sigma\in\Sigma}AC(E^{\sigma})^{*}$
to $\mathcal{B}:=\coprod_{\sigma\in\Sigma}B(H_{\sigma})$ with the
property that $\mathbf{\widehat{X}}$ maps $AC(E^{\sigma})^{*}$ to
$B(H_{\sigma})$. In fact, it will be convenient to view $\mathcal{B}$
as a bundle over $\mathcal{AC}(E)$ with the property that the total
space of $\mathcal{B}|_{AC(E^{\sigma})^{*}}$ is $AC(E^{\sigma})^{*}\times B(H_{\sigma})$.
When this is done, we follow the customary practice of identifying
a section of a trivial bundle over a space with a function on the
space with values in the fibre. (In this case the trivial bundle is
$AC(E^{\sigma})^{*}\times B(H_{\sigma})$ over $AC(E^{\sigma})^{*}$.)
Then we can say, simply, that $\mathbf{\widehat{X}}$ is a \emph{section}
of this bundle and adopt the following terminology.

\begin{definition}The section $\mathbf{\widehat{X}}$ of the bundle
$\mathcal{B}=\coprod_{\sigma\in\Sigma}AC(E^{\sigma})^{*}\times B(H_{\sigma})$
over $\mathcal{AC}(E)=\coprod_{\sigma\in\Sigma}AC(E^{\sigma})^{*}$
associated with the element $X\in H^{\infty}(E)$ is called \emph{the
complete Schur class section} determined by $X$. \end{definition}

\begin{remark}It is a consequence of \cite[Lemma 3.8]{MSSchur} and
the fact that we incorporate \emph{all} the normal representations
of $M$ that the map $X\to\mathbf{\widehat{X}}$ is injective. It
is very much of interest to understand how to adjust matters when
one restricts attention to some subset of $\Sigma$. \end{remark}

A natural question is, ``How does one recognize such a section?''.
The answer begins with the structure of the families of sets $\mathcal{AC}(E)=\{AC(E^{\sigma}):\sigma\in\Sigma\}$
and $\mathcal{D}(E):=\{\mathbb{D}(E^{\sigma}):\sigma\in\Sigma\}$,
which is abstracted by the following definition.

\begin{definition} A family $\mathcal{A}=\{\mathcal{A}(\sigma):\sigma\in\Sigma\}$
is said to be a \emph{fully matricial} \emph{E-set} if
\begin{enumerate}
\item [(i)] for each $\sigma$, $\mathcal{A}(\sigma)\subseteq E^{\sigma}$,
\item [(ii)] it is closed with respect to taking direct sums; that is,
$\mathcal{A}(\sigma)\oplus\mathcal{A}(\tau)\subseteq\mathcal{A}(\sigma\oplus\tau)$
and
\item [(iii)] it is closed with respect to unitary similarity; that is,
if $\eta\in\mathcal{A}(\sigma)$ and $u\in\sigma(M)'$ is a unitary
then $u\cdot\eta\cdot u^{*}\in\mathcal{A}(\sigma)$.
\end{enumerate}
\end{definition}

Note that the product $u\cdot\eta\cdot u^{*}$ is the bimodule product
on $E^{\sigma}$.

Our notion is very similar to one defined by Taylor in \cite{jT72}
and to one defined by Voiculescu in \cite{Voi}. In these papers,
the general linear group is used instead of the unitary group. For
our purposes here, however, it is more convenient to work with the
unitary group, as Helton, Klepp and McCullough did in \cite{HKM}.
It is evident that just as with $\mathcal{A}C(E)$, we may view $\mathcal{B}$
as a bundle over any fully matricial $E$-set and study sections of
this bundle. The following definition singles out a special property
that a section may or may not have.

First, we need to extend the definition of the spaces $\mathcal{I}(S_{0},\eta^{*})$
that we have met before: For \emph{any} two $\sigma,\tau\in\Sigma$,
and $\eta\in E^{\sigma}$ and $\zeta\in E^{\tau}$, we write $\mathcal{I}(\eta^{*},\zeta^{*})$
for the collection of all $C:H_{\sigma}\rightarrow H_{\tau}$ such
that $C\sigma(\cdot)=\tau(\cdot)C$ and such that $C\eta^{*}=\zeta^{*}(I_{E}\otimes C)$,
as maps from $E\otimes_{\sigma}H_{\sigma}$ to $H_{\tau}$. We call
an element of $\mathcal{I}(\eta^{*},\zeta^{*})$ an intertwiner of
$\eta^{*}$ and $\zeta^{*}$. It is easy to see that if $\eta$ and
$\zeta$ both have norm at most one, then $\mathcal{I}(\eta^{*},\zeta^{*})$
is simply the collection of operators that intertwine $\sigma\times\eta^{*}$
and $\tau\times\zeta^{*}$. Note that an intertwiner $C\in\mathcal{I}(\eta^{*},\zeta^{*})$
will also satisfies $C\eta_{k}^{*}=\zeta_{k}^{*}(I_{E^{\otimes k}}\otimes C)$
for all $k\geq1$ where, recall, $\eta_{k}^{*}$ and $\zeta_{k}^{*}$
are the generalized powers of $\eta^{*}$ and $\zeta^{*}$ discussed
above.

\begin{definition}\label{freehol}

Let $\mathcal{A}=\{\mathcal{A}(\sigma):\sigma\in\Sigma\}$ be a fully
matricial $E$-set and form the bundle $\mathcal{B}=\coprod_{\sigma\in\Sigma}\mathcal{A}(\sigma)\times B(H_{\sigma})$.
We say that a section $\mathbf{f}$ of $\mathcal{B}$ \emph{preserves
intertwiners} in case $Cf_{\sigma}(\eta)=f_{\tau}(\zeta)C$ for all
$\sigma,\tau\in\Sigma$, $(\eta,\zeta)\in\mathcal{A}(\sigma)\times\mathcal{A}(\tau)$,
and intertwiners $C\in\mathcal{I}(\eta^{*},\zeta^{*})$.

\end{definition} With all the pieces before us, we may formulate
our ``recognition'' theorem as follows.

\begin{theorem}\label{Schur_v_Hinfty} A section $\mathbf{f}$ of
the bundle $\mathcal{B}$ over $\mathcal{AC}(E)$ is a complete Schur
section if and only if $\mathbf{f}$ preserves intertwiners. \end{theorem}

We won't go into the details of the proof here, but we do want to
point out that a special role is played by the fact that the functions
are defined on $\mathcal{AC}$ and not just on the the family of open
sets $\mathcal{D}(E)=\coprod_{\sigma\in\Sigma}\mathbb{D}(E^{\sigma})^{*}$.
A key role is played by our Theorem \ref{AbsContConditions}; the
interaction between ``absolute continuity'' and the special nature
of elements in $H^{\infty}(E)$ is what underlies the proof.

Now a section $\mathbf{\widehat{X}}$, $X\in H^{\infty}(E)$, may
be restricted to $\mathcal{D}(E)$ and when this is done, the resulting
section is analytic as a Banach-space-valued section, i.e., $\widehat{X_{\sigma}}$
is a $B(H_{\sigma})$-valued analytic function on $\mathbb{D}(E^{\sigma})^{*}$.
But there are many sections of $\mathcal{B}|_{\mathcal{D}(E)}$ with
this property that don't come from elements of $H^{\infty}(E)$. Here
is a very simple example.

{}

\begin{example} Let $M=E=\mathbb{C}$. In this case $H^{\infty}(E)$
is the classical $H^{\infty}(\mathbb{D})$. The representations in
$\Sigma$ are just the obvious representations. We let $\sigma$ be
the identity representation of $\mathbb{C}$ on $\mathbb{C}$. Then
every representation of $\mathbb{C}$ is a multiple of $\sigma$,
$n\sigma$, which acts on $\mathbb{C}^{n}$. We treat $\mathbb{C}^{\infty}$
as $\ell^{2}(\mathbb{N})$. Then $E^{n\sigma}=E^{n\sigma*}=B(\mathbb{C}^{n})$
and $\mathbb{D}(E^{n\sigma})^{*}=\{A\in B(\mathbb{C}^{n})\mid\norm{A}<1\}$.
We set $f_{n\sigma}(A)=(I-A)^{-1}$, for $A\in\mathbb{D}(E^{n\sigma})^{*}$.

If $A\in B(\mathbb{C}^{n})$, $B\in B(\mathbb{C}^{m})$ both have
norm less than $1$ and if $C:\mathbb{C}^{m}\rightarrow\mathbb{C}^{n}$
intertwines them, that is, if $AC=CB$, then $C$ also intertwines
$(I-A)^{-1}$ and $(I-B)^{-1}$. Thus $\mathbf{f}=\{f_{n\sigma}\}$
is a section of $\mathcal{B}|_{\mathcal{D}(E)}$ which certainly deserves
to be called analytic. After all, it comes from the function $h$,
where $h(z)=\sum_{n\geq0}z^{n}$. However, $h$ is not in $H^{\infty}(\mathbb{D})$.

\end{example}

Note that in this example, we have identified $H^{\infty}(\mathbb{D})$,
which is a space of complex valued functions on the (classical) unit
disc $\mathbb{D}=\mathbb{D}(E^{\sigma})^{*}$, with $H^{\infty}(\mathbb{C})$,
which really is a space of sequences, viz., the space of the sequences
of Taylor coefficients of the functions in $H^{\infty}(\mathbb{D})$.
In the case of $h$, of course, the sequence is $(1,1,1,\cdots)$,
which is not in $H^{\infty}(\mathbb{C})$. In general, recall, every
element $X\in H^{\infty}(E)$ has a series expansion \eqref{eq:Taylor series.},
so it is natural to wonder if it is possible to manipulate arbitrary
series of tensors. It \emph{is} possible, and to help clarify how,
we introduce the following definition.

\begin{definition}Let $E$ be a $W^{*}$-correspondence over the
$W^{*}$-algebra $M$.
\begin{enumerate}
\item A \emph{(formal) series of tensors} (over $E)$ is simply a sequence
$\theta=\{\theta_{k}\}_{k\ge0}$, where $\theta_{k}\in E^{\otimes k}$.
However, we shall usually write $\theta\sim\sum_{k\geq0}\theta_{k}$
in anticipation of function-theoretic considerations to come.
\item If $\theta\sim\sum_{k\geq0}\theta_{k}$ is a series of tensors over
$E$, then we define $R(\theta)$ to be \[
(\overline{\lim}_{k}\norm{\theta_{k}}^{1/k})^{-1},\]
and we refer to $R(\theta)$ as \emph{the radius of convergence} of
$\theta$.
\end{enumerate}
\end{definition}

Evidently, $R(\theta)$ is a non-negative number or $+\infty$. The
formula for $R(\theta)$ suggests that $\Sigma_{k\geq0}\theta_{k}$
converges in some sense. And of course it does, as the following theorem
shows. It is a generalization of the well known Cauchy-Hadamard theorem
from elementary complex analysis. It plays a prominent role in Popescu's
study of free analyticity, also.

\begin{theorem}\label{RadiusConv} Suppose $\theta\sim\sum_{k\geq0}\theta_{k}$
is a series of tensors coming from $E$, and let $R=R(\theta)$ be
its radius of convergence.
\begin{enumerate}
\item [(1)] Given $\sigma\in\Sigma$ and $\eta^{*}\in R\mathbb{D}(E^{\sigma})^{*}:=\{R\zeta^{*}\mid\zeta^{*}\in\mathbb{D}(E^{\sigma})^{*}\}$,
the series $\sum_{k}\norm{\eta_{k}^{*}L_{\theta_{k}}}$ converges,
where, recall, $\eta_{k}^{*}$ denotes the $k^{th}$ generalized power
of $\eta^{*}$ and where $L_{\theta_{k}}$is the map from $H$ to
$E^{\otimes k}\otimes_{\sigma}H$ defined by $L_{\theta_{k}}h=\theta_{k}\otimes h$.
If $0<\rho<R$, the convergence is uniform on $\rho\overline{\mathbb{D}(E^{\sigma})}$.
\item [(2)] If $R<R'<\infty$, there exists a $\sigma\in\Sigma$ and an
$\eta\in E^{\sigma}$ with $\norm{\eta}=R'$ such that $\sum_{k}\norm{\eta_{k}^{*}L_{\theta_{k}}}=\infty$.
\end{enumerate}
\end{theorem}

Thus, a series $\theta$ defines a function on each of the discs $R\mathbb{D}(E^{\sigma})^{*}$
with values values in $B(H_{\sigma})$. We denote this function by
$\widehat{\theta_{\sigma}}$. Its value at an $\eta^{*}$ is given
by the formula,\[
\widehat{\theta_{\sigma}}(\eta^{*})=\sum_{k\geq0}\eta_{k}^{*}L_{\theta_{k}}.\]
By Theorem \ref{RadiusConv}, this series is a series of operators
$B(H_{\sigma})$ that converges in norm. The family $\{\widehat{\theta_{\sigma}}\}_{\sigma\in\Sigma}$
forms the section $\mathbf{\widehat{\theta}}$ of the bundle $\mathcal{B}:=\coprod_{\sigma\in\Sigma}R\mathbb{D}(E^{\sigma})^{*}\times B(H_{\sigma})$
over $R\mathcal{D}(E):=\coprod_{\sigma\in\Sigma}R\mathbb{D}(E^{\sigma})^{*}$.

\begin{definition}Let $\theta\sim\sum_{k\geq0}\theta_{k}$ be a series
of tensors over $E$ and let $R$ be its radius of convergence. The
section $\mathbf{\widehat{\theta}}$ of the bundle $\mathcal{B}:=\coprod_{\sigma\in\Sigma}R\mathbb{D}(E^{\sigma})^{*}\times B(H_{\sigma})$
over $R\mathcal{D}(E):=\coprod_{\sigma\in\Sigma}R\mathbb{D}(E^{\sigma})^{*}$determined
by the family $\{\widehat{\theta_{\sigma}}\}_{\sigma\in\Sigma}$ is
called \emph{the free analytic section} determined by $\theta$. \end{definition}

We may now characterize the free analytic sections in much the same
fashion as we characterized complete Schur sections.

\begin{theorem}\label{RD} A section $\mathbf{f}=\{f_{\sigma}\}_{\sigma\in\Sigma}$
of the bundle
\[\mathcal{B}=\coprod_{\sigma\in\Sigma}R\mathbb{D}(E^{\sigma})^{*}\times
B(H_{\sigma}) \mbox{ over } R\mathcal{D}(E)\] is the free analytic section
determined by
a series of tensors with radius of convergence at least $R$ if and
only if \textbf{$\mathbf{f}$} preserves intertwiners. \end{theorem}


\end{document}